\documentclass[11pt]{amsart}

\usepackage{amsmath}
\usepackage{amscd}
\usepackage{amssymb}
\usepackage{latexsym}
\usepackage{stmaryrd} 
\usepackage{color}
\usepackage{mathrsfs}

\usepackage[arrow,curve,matrix,tips,2cell]{xy}
  \SelectTips{eu}{10} \UseTips
  \UseAllTwocells
\newtheorem{theorem}{Theorem}[section]
\newtheorem*{theorem*}{Theorem}
\newtheorem{lemma}[theorem]{Lemma}
\newtheorem*{lemma*}{Lemma}
\newtheorem{corollary}[theorem]{Corollary}
\newtheorem*{corollary*}{Corollary}
\newtheorem{proposition}[theorem]{Proposition}

\newtheorem{remark}[theorem]{Remark}
\newtheorem{definition}[theorem]{Definition}

\newtheorem{question}[theorem]{Question}
\newcommand{\bgl}{\begin{equation}} 
\newcommand{\egl}{\end{equation}}
\newcommand{\bgloz}{\begin{equation*}} 
\newcommand{\egloz}{\end{equation*}}
\newcommand{\bgln}{\begin{eqnarray}} 
\newcommand{\egln}{\end{eqnarray}}
\newcommand{\bglnoz}{\begin{eqnarray*}} 
\newcommand{\eglnoz}{\end{eqnarray*}}
\newcommand{\btheo}{\begin{theorem}}
\newcommand{\etheo}{\end{theorem}}
\newcommand{\btheooz}{\begin{theorem*}}
\newcommand{\etheooz}{\end{theorem*}}
\newcommand{\blemma}{\begin{lemma}}
\newcommand{\elemma}{\end{lemma}}
\newcommand{\blemmaoz}{\begin{lemma*}}
\newcommand{\elemmaoz}{\end{lemma*}}
\newcommand{\bproof}{\begin{proof}}
\newcommand{\eproof}{\end{proof}}
\newcommand{\bbew}{\begin{beweis}}
\newcommand{\ebew}{\end{beweis}}
\newcommand{\bremark}{\begin{remark}\em}
\newcommand{\eremark}{\end{remark}}
\newcommand{\bdefin}{\begin{definition}}
\newcommand{\edefin}{\end{definition}}
\newcommand{\bprop}{\begin{proposition}}
\newcommand{\eprop}{\end{proposition}}
\newcommand{\bcor}{\begin{corollary}}
\newcommand{\ecor}{\end{corollary}}
\newcommand{\bcoroz}{\begin{corollary*}}
\newcommand{\ecoroz}{\end{corollary*}}
\newcommand{\bfa}{\begin{cases}} 
\newcommand{\efa}{\end{cases}}
\newcommand{\bquestion}{\begin{question}\em}
\newcommand{\equestion}{\end{question}}
\newcommand{\cE}{\mathcal E}
\newcommand{\cF}{\mathcal F}

\newcommand{\cI}{\mathcal I}
\newcommand{\cJ}{\mathcal J}

\newcommand{\cL}{\mathcal L}

\newcommand{\cP}{\mathcal P}
\newcommand{\cQ}{\mathcal Q}
\newcommand{\cR}{\mathcal R}

\newcommand{\cV}{\mathcal V}

\def\Az{\mathbb{A}}
\def\Cz{\mathbb{C}}
\def\Fz{\mathbb{F}}

\def\Nz{\mathbb{N}}

\def\Qz{\mathbb{Q}}
\def\Rz{\mathbb{R}}

\def\Zz{\mathbb{Z}}

\def\1z{\mathbb{1}}
\newcommand{\fJ}{\mathfrak J}

\newcommand{\mfa}{\mathfrak a}

\newcommand{\mfp}{\mathfrak p}
\newcommand{\sR}{\mathscr{R}}

\newcommand{\an}[1]{``#1''} 
\newcommand{\ti}{\tilde}
\newcommand{\wt}{\widetilde}

\newcommand{\lori}{\longrightarrow}

\newcommand{\ma}{\mapsto} 
\newcommand\into{\hookrightarrow} 
\newcommand{\LRarr}{\Leftrightarrow} 

\newcommand{\ve}{\varepsilon}

\def\SEMI{\mbox{$\times\kern-2pt\vrule height5pt width.6pt \kern3pt $}}

\newcommand{\halb}{\tfrac{1}{2}}

\newcommand{\ind}{{\rm ind\,}}
\newcommand{\res}{{\rm res\,}}

\newcommand{\id}{{\rm id}}

\renewcommand{\ker}{{\rm ker}\,}
\newcommand{\coker}{{\rm coker}\,}
\newcommand{\im}{{\rm im}\,}
\newcommand{\reg}{^\times} 
\newcommand{\pos}{_{>0}} 
\newcommand{\lspan}{{\rm span}} 
\newcommand{\abs}[1]{\lvert#1\rvert} 
\newcommand{\defeq}{\mathrel{:=}} 
\newcommand{\eqdef}{\mathrel{=:}} 

\newcommand{\dop}{\text{: }} 
\newcommand{\ilim}{\varinjlim} 

\newcommand{\Proj}{{\rm Proj}\,}
\newcommand{\Env}{{\rm Env}\,}

\newcommand{\tight}{{\rm tight}}

\newcommand{\lge}{\left\{} 
\newcommand{\rge}{\right\}} 
\newcommand{\lru}{\left(} 
\newcommand{\rru}{\right)} 
\newcommand{\leck}{\left[} 
\newcommand{\lsp}{\left\langle} 
\newcommand{\rsp}{\right\rangle} 
\newcommand{\rukl}[1]{\lru #1 \rru} 
\newcommand{\gekl}[1]{\lge #1 \rge} 
\newcommand{\spkl}[1]{\lsp #1 \rsp} 
\newcommand{\menge}[2]{\gekl{ #1 \dop #2 }} 
\begin{document}

\title[Independent resolutions II: C*-algebraic case]{Independent resolutions for totally disconnected dynamical systems II: C*-algebraic case}

\author{Xin Li}
\author{Magnus Dahler Norling}

\address{Xin Li, School of Mathematical Sciences, Queen Mary University of London, Mile End Road, London E1 4NS, UK}
\email{xin.li@qmul.ac.uk}
\address{Magnus Dahler Norling, Institute of Mathematics, University of Oslo, P.b. 1053 Blindern, 0316 Oslo, Norway}
\email{magnudn@math.uio.no}


\thanks{\scriptsize{The research of the first named author was partially supported by the ERC through AdG 267079.}}

\begin{abstract}
We develop the notion of independent resolutions for crossed products attached to totally disconnected dynamical systems. If such a crossed product admits an independent resolution of finite length, then its K-theory can be computed (at least in principle) by analysing the corresponding six-term exact sequences. Building on our previous paper on algebraic independent resolutions, we give a criterion for the existence of finite length independent resolutions. Moreover, we illustrate our ideas in various concrete examples.
\end{abstract}

\maketitle


\setlength{\parindent}{0pt} \setlength{\parskip}{0.5cm}

\section{Introduction}

The crossed product construction is one of the most classical constructions in operator algebras, and topological K-theory is one of the most important invariants for C*-algebras. Therefore, a very natural task is to find systematic ways to compute K-theory for C*-algebraic crossed products.

The goal of the present paper is to take up this task in the situation of crossed products attached to totally disconnected dynamical systems. We do so using the central notion of independent resolutions. In our previous paper \cite{L-N}, we introduced and discussed independent resolutions from a purely algebraic point of view. Now, our goal is to develop a notion of independent resolutions in the C*-algebraic setting. Building on our previous work \cite{L-N}, we then produce C*-algebraic independent resolutions which allow us to compute K-theory for crossed products. More precisely, let $\Omega$ be a totally disconnected locally compact Hausdorff space and $\Gamma$ a discrete group acting on $\Omega$. Consider the reduced crossed product $C_0(\Omega) \rtimes_r \Gamma$. If $\Gamma$ satisfies the Baum-Connes conjecture with coefficients and $\Omega$ admits a $\Gamma$-invariant regular basis in the sense of \cite{C-E-L2} (see also \S~\ref{sec-ind-res} for explanations), then the main result in \cite{C-E-L2} provides a formula for the K-theory of $C_0(\Omega) \rtimes_r \Gamma$. However, it was also observed in \cite{C-E-L2} that in general, it is not possible to find a $\Gamma$-invariant regular basis. Still, following \cite[Remark~3.22]{C-E-L2}, what we can always do is to produce a sequence $X$, $X_1$, $X_2$, ... of totally disconnected $\Gamma$-spaces which admit $\Gamma$-invariant regular bases and which fit into a $\Gamma$-equivariant long exact sequence $\dotso \to C_0(X_2) \to C_0(X_1) \to C_0(X) \to C_0(\Omega) \to 0$. We call this an independent resolution of $\Gamma \curvearrowright C_0(\Omega)$. Under the assumption that $\Gamma$ is exact, the sequence $\dotso \to C_0(X_2) \rtimes_r \Gamma \to C_0(X_1) \rtimes_r \Gamma \to C_0(X) \rtimes_r \Gamma \to C_0(\Omega) \rtimes_r \Gamma \to 0$ will still be exact, and we call this an independent resolution of $C_0(X) \rtimes_r \Gamma$. If, furthermore, $\Gamma$ satisfies the Baum-Connes conjecture with coefficients, then we can apply the K-theoretic formula from \cite{C-E-L2} to each of the crossed products $C_0(X) \rtimes_r \Gamma$, $C_0(X_1) \rtimes_r \Gamma$, ... and try to compute K-theory for $C_0(\Omega) \rtimes_r \Gamma$ using our long exact sequence. In general, given an independent resolution of $\Gamma \curvearrowright C_0(\Omega)$ satisfying a certain freeness condition for the group actions, there is at least a spectral sequence which converges to $K_*(C_0(\Omega) \rtimes_r \Gamma)$ in good cases.

The case where we have a finite length independent resultion (i.e., we can choose $X_{n+1} = \emptyset$ for some $n$) is particularly nice. In that case, the exact sequence $0 \to C_0(X_n) \rtimes_r \Gamma \to \dotso \to C_0(X_2) \rtimes_r \Gamma \to C_0(X_1) \rtimes_r \Gamma \to C_0(X) \rtimes_r \Gamma \to C_0(\Omega) \rtimes_r \Gamma \to 0$ splits into short exact sequences which can be studied in K-theory by means of six-term exact sequences. The point is that given a finite length independent resolution, we only have to solve finitely many six-term exact sequences. And if we try to solve these successively, we will always be in the situation that we already know the K-groups for two out of the three C*-algebras which appear in each of our sequences.

The main goal of this paper is to give a criterion which guarantees the existence of finite length independent resolutions. This builds on \cite{L-N}. The bridge between algebraic independent resolutions and C*-algebraic ones is given by the observation that a sequence of totally disconnected dynamical systems which all admit invariant regular bases gives rise to an algebraic independent resolution if and only if it gives rise to a C*-algebraic one. In addition, these independent resolutions are intimately related. For instance, the homomorphisms in the algebraic independent resolution induce the ones in the C*-algebraic independent resolution. In particular, the former one has finite length if and only if the latter one does. Therefore, the criterion for the existence of finite length algebraic independent resolutions in \cite{L-N} gives us a criterion for the existence of C*-algebraic independent resolutions of finite length.

We remark that finding such an independent resolution of finite length is only the first step in the K-theory computation for our crossed product. The second step is to go through the short exact sequences into which our exact sequence splits and to compute all the corresponding six-term exact sequences. It might be that we encounter serious extension problems along the way, so this second step might require extra work.

In order to illustrate our main result, we discuss various concrete examples. If we want to apply our ideas to compute K-theory for a given C*-algebra, the first step is to describe the C*-algebra as a crossed product of a totally disconnected dynamical system, at least up to Morita equivalence. This is for instance possible for C*-algebras of certain $0$-$F$ inverse semigroups and certain quotients of these. This has already been observed in \cite{Nor2}, but we present a slightly different approach which is more explicit and better suited for our purposes. More concretely, we discuss graph C*-algebras and one dimensional tiling C*-algebras, and derive crossed product descriptions for these. This might be of independent interest. We then use independent resolutions to compute K-theory for graph C*-algebras and C*-algebras of one dimensional tilings. We also determine K-theory for certain ideals and quotients of semigroup C*-algebras. In particular, our method allows us to study the K-theory of group C*-algebras with the help of semigroup C*-algebras. The idea is to choose a suitable subsemigroup of our group which gives rise to a finite length independent resolution for the group C*-algebra we are interested in. Furthermore, our ideas allow us to compute K-theory for the C*-algebra of semigroups which do not satisfy the independence condition. Such semigroups could not be treated using the original method of \cite{C-E-L2}. Interestingly, in our example, we again encounter the phenomenon that the K-theories of the left and right reduced semigroup C*-algebras coincide.

\section{Independent resolutions}
\label{sec-ind-res}

The notion of independent resolutions has already been introduced in \cite{L-N}, but in a purely algebraic setting. We now discuss C*-algebraic independent resolutions.

Throughout this paper, every group is supposed to be discrete and countable, and every topological space is assumed to be second countable, locally compact and Hausdorff. Given a dynamical system $\Gamma \curvearrowright \Omega$ with a group $\Gamma$ acting on a totally disconnected space $\Omega$, we want to introduce the notion of an independent resolution of $\Gamma \curvearrowright C_0(\Omega)$. Once we have done that, we can also talk about independent resolutions for dynamical systems of the form $\Gamma \curvearrowright D$ where $D$ is a commutative C*-algebra generated by projections since C*-algebras of the form $C_0(\Omega)$ for a totally disconnected space $\Omega$ are precisely those commutative C*-algebras which are generated by projections.

First of all, a semilattice is by definition a commutative idempotent semigroup, i.e., a commutative semigroup in which every element $e$ satisfies $e e = e$. All our semilattices are supposed to have a zero element. Given a semilattice $E$, the C*-algebra of $E$ is the universal C*-algebra
$$ 
C^*_u(E) = C^* \rukl{ \gekl{p_e}_{e \in E} \vline 
  \begin{array}{c}
  p_e \text{ are projections}, \ p_0 = 0,\\
  E \ni e \ma p_e \text{ is a semigroup homomorphism}
  \end{array}
  }
$$
By an action of a group $\Gamma$ on a semilattice $E$ we mean a group homomorphism from $\Gamma$ to the semigroup automorphisms of $E$. Such an action obviously induces an action of $\Gamma$ on $C^*_u(E)$.

It turns out that every C*-algebra of the form $C_0(\Omega)$ for a totally disconnected space $\Omega$ is isomorphic to the C*-algebra of a suitable semilattice. Namely, by \cite[Proposition~2.12]{C-E-L2}, we can always find a regular basis $\cV$ for $\Omega$ in the sense of \cite[Definition~2.9]{C-E-L2}. Since the compact open sets in $\cV$ are closed under intersection, they form a semilattice. And as explained in \cite[Remark~3.22]{C-E-L2}, we have the isomorphism $C^*_u(\cV) \cong C_0(\Omega)$, $p_V \ma 1_V$. Here $p_V$ is the projection in the C*-algebra of our semilattice $\cV$ corresponding to $V \in \cV$ (as in the definition of $C^*_u(\cV)$), and $1_V$ is the characteristic function of $V$.

Now given a totally disconnected dynamical system $\Gamma \curvearrowright \Omega$, we can ask for a semilattice $E$, together with an action of $\Gamma$, such that we have a $\Gamma$-equivariant isomorphism $C^*_u(E) \cong C_0(\Omega)$. It is easy to see that such a system $\Gamma \curvearrowright E$ exists for $\Gamma \curvearrowright \Omega$ if and only if $\Omega$ admits a $\Gamma$-invariant regular basis in the sense of \cite[Definition~2.9]{C-E-L2}. In general, this does not need to be the case, as was remarked in \cite[Proposition~3.18]{C-E-L2}. However, \cite[Remark~3.22]{C-E-L2} shows that given an arbitrary totally disconnected dynamical system $\Gamma \curvearrowright \Omega$, we can always find semilattices $E$, $E_1$, $E_2$, ..., together with $\Gamma$-actions on these semilattices, and a $\Gamma$-equivariant long exact sequence
\bgl
\label{ind-res}
  \dotso \to C^*_u(E_2) \to C^*_u(E_1) \to C^*_u(E) \to C_0(\Omega) \to 0.
\egl
We call such a long exact sequence an independent resolution of $\Gamma \curvearrowright C_0(\Omega)$. Of course, the requirement that the sequence is $\Gamma$-equivariant is crucial here. Moreover, we define the length of such an independent resolution to be the smallest integer $n \geq 0$ with $E_{n+1} = \gekl{0}$, or equivalently, $C^*_u(E_{n+1}) = \gekl{0}$. If no such integer exists, then we set the length to be $\infty$.

An independent resolution of $\Gamma \curvearrowright C_0(\Omega)$ for which the stabilizer groups are all trivial ($\Gamma$ acts freely on $E\reg$ and $E_k\reg$ for all $k$) is a $\fJ$-projective resolution for $C_0(\Omega)$ in the category $KK^{\Gamma}$, in the sense of \cite[\S~2]{Mey}. Here we take for $\fJ$ the K-theory functor from the category $KK^{\Gamma}$ to $\Zz / 2 \Zz$-graded $\Zz \Gamma$-modules. As explained in \cite[\S~3]{Mey}, every $\fJ$-projective resolution embeds into a phantom tower, which in turn induces the so-called ABC spectral sequence \cite[\S~4]{Mey}. In \cite{Mey}, the reader may find conditions under which this ABC spectral sequence converges to $K_*(C_0(\Omega) \rtimes_r \Gamma)$ (see for instance \cite[Proposition~4.1]{Mey} or \cite[\S~5]{Mey}). The reader may find more details in \cite{M-N} and \cite{Mey}. But at least in principle, an independent resolution with trivial stabilizer groups helps to compute the K-theory of our crossed product. In the case of finite length resolutions, we elaborate on this computational aspect in \S~\ref{k-theory-free-action-section}.

Now let us assume that our group $\Gamma$ is exact. In that case, every independent resolution as in \eqref{ind-res} gives rise to a long exact sequence of the form
\bgl
\label{les}
  \dotso \to C^*_u(E_2) \rtimes_r \Gamma \to C^*_u(E_1) \rtimes_r \Gamma \to C^*_u(E) \rtimes_r \Gamma \to C_0(\Omega) \rtimes_r \Gamma \to 0.
\egl
Here we take the crossed products with respect to the $\Gamma$-actions provided by our independent resolution. We call such a long exact sequence an independent resolution of $C_0(\Omega) \rtimes_r \Gamma$. As remarked at the beginning, we can also talk about independent resolutions for dynamical systems of the form $\Gamma \curvearrowright D$ or for $D \rtimes_r \Gamma$ where $D$ is a commutative C*-algebra generated by projections.

If $\Gamma \curvearrowright C_0(\Omega)$ admits an independent resolution of finite length, then we get the following exact sequence:
$$
  0 \to C^*_u(E_n) \rtimes_r \Gamma \to \dotso \to C^*_u(E_1) \rtimes_r \Gamma \to C^*_u(E) \rtimes_r \Gamma \to C_0(\Omega) \rtimes_r \Gamma \to 0.
$$
This exact sequence can be split into several short exact sequences of the form
\bglnoz
  && 0 \to C^*_u(E_n) \rtimes_r \Gamma \to C^*_u(E_{n-1}) \rtimes_r \Gamma \to \mathrm{ker}_{n-2} \to 0 \\
  && 0 \to \mathrm{ker}_{n-2} \to C^*_u(E_{n-2}) \rtimes_r \Gamma \to \mathrm{ker}_{n-3} \to 0 \\
  && \dots \\
  && 0 \to \mathrm{ker}_1 \to C^*_u(E_1) \rtimes_r \Gamma \to \mathrm{ker}_0 \to 0 \\
  && 0 \to \mathrm{ker}_0 \to C^*_u(E) \rtimes_r \Gamma \to C_0(\Omega) \rtimes_r \Gamma \to 0.
\eglnoz
Now consider the corresponding six-term exact sequences in K-theory, and assume that $\Gamma$ satisfies the Baum-Connes conjecture with coefficients. In the first six-term exact sequence, the K-theories for $C^*_u(E_n) \rtimes_r \Gamma$ and $C^*_u(E_{n-1}) \rtimes_r \Gamma$ can be computed using \cite[Corollary~3.14]{C-E-L2}. If it is possible to compute the K-theory for $\mathrm{ker}_{n-2}$ from this six-term exact sequence, we could plug in the result into the next six-term exact sequence, apply \cite[Corollary~3.14]{C-E-L2} to $C^*_u(E_{n-2}) \rtimes_r \Gamma$, and try to determine the K-theory of $\mathrm{ker}_{n-3}$. In this way, we could compute K-theory step by step until we come to the C*-algebra of interest, namely $C_0(\Omega) \rtimes_r \Gamma$. Of course, the extension problems which we have to solve along the way might be difficult.

\section{From algebraic independent resolutions to independent resolutions}
\label{alg-Cstalg}

Let us now build the bridge between algebraic independent resolutions and C*-algebraic ones.

Let $A$ be a C*-algebra generated by a multiplicatively closed family of projections $\cP$, and let $Z$ be the sub-$\Zz$-algebra of $A$ generated by $\cP$. Assume that $E$ is a semilattice with a semilattice homomorphism $E \to \cP$, which induces homomorphisms $\pi_\Zz$: $\Zz_0[E] \to Z$ and $\pi$: $C^*_u(E) \to A$. Let $E'$ be a semilattice of projections in $\Zz_0[E]$, let $I_{\Zz} = \Zz \text{-} \lspan(E')$ and $I$ be the ideal of $C^*_u(E)$ generated by $I_{\Zz}$.

\blemma
If $\ker \pi_{\Zz} = I_{\Zz}$, then $\ker \pi = I$.
\elemma
\bproof
Let $\cF$ be the collection of finite subsets of $E'$ which are closed under multiplication. $\cF$ is obviously inductively ordered with respect to inclusion. Moreover, set for $F \in \cF$: $C^*_F(E) \defeq C^*(\menge{e}{e \in F}) \subseteq C^*_u(E)$. We obviously have $C^*_u(E) = \overline{\bigcup_{F \in \cF} C^*_F(E)}$. Since $C^*_u(E) / I = \overline{\bigcup_{F \in \cF} (C^*_F(E) / I_F)}$ with $I_F = C^*_F(E) \cap I$, all we have to prove is that the homomorphism induced by restricting $\pi$ to $C^*_F(E)$, $\pi \vert_F : \ C^*_F(E) / I_F \to A$, is injective for all $F \in \cF$. Given $F \in \cF$, we can orthogonalize the projections in $F$ and obtain a new set of non-zero projections $F^{(\text{orth})}$. But since $F$ is multiplicatively closed, we have $F^{(\text{orth})} \subseteq \Zz \text{-} \lspan (F)$. Since $C^*_F(E) = \bigoplus_{f \in F^{(\text{orth})}} \Cz \cdot f$, $\pi \vert_F$ is injective if and only if for all $f \in F^{(\text{orth})}$, $\pi(f) = 0$ implies $f \in I_F$. But $\pi(f) = 0$ means that $\pi_{\Zz}(f) = 0$, so that $f \in \ker \pi_{\Zz}$. By assumption, $f$ must lie in $I_{\Zz}$. Hence $f \in I \cap C^*_F(E) = I_F$, as desired.
\eproof

\bcor
\label{algindres-->Cstarindres}
Let $\Gamma \curvearrowright \Omega$ be a totally disconnected dynamical system. Assume that $E$, $E_1$, $E_2$, ... are $\Gamma$-semilattices and that
\bgloz
  \dotso \to \Zz_0[E_2] \to \Zz_0[E_1] \to \Zz_0[E] \to C_0(\Omega, \Zz) \to 0
\egloz
is an algebraic independent resolution. Then
\bgloz
  \dotso \to C^*_u(E_2) \to C^*_u(E_1) \to C^*_u(E) \to C_0(\Omega) \to 0
\egloz
is an independent resolution. The homomorphisms in this sequence are induced by the ones from the algebraic independent resolution.
\ecor

In the following, we give a criterion for the existence of finite length independent resolutions. The previous corollary reduces our investigations to the algebaic setting, so that we can use \cite[\S~4]{L-N}. We first introduce some notation. Let $E$ be a semilattice. A finite cover for $e \in E\reg$ is a finite subset $\gekl{f_j}_{j \in J}$ of $E\reg$ ($J$ is a finite index set) with the property that
\begin{itemize}
\item $f_j \leq e$ for all $j \in J$,
\item for every $f \in E\reg$ with $f \leq e$, there exists $j \in J$ such that $f f_j \neq 0$.
\end{itemize}
Given a finite cover $\gekl{f_j}_{j \in J}$ for $e \in E\reg$, we can write, in a unique way, $\bigvee_{j \in J} f_j = \sum_{k} n_k \ve_k$ where the $\ve_k$ are pairwise distinct idempotents in $E\reg$ and the $n_k$ are non-zero integers. Here $\bigvee_{j \in J} f_j$ is the smallest projection in $\Zz_0[E]$ which dominates all the $f_j$. We set $\bigvee\gekl{f_j}_{j \in J} \defeq \bigvee_{j \in J} f_j$ and $E(\bigvee\gekl{f_j}_{j \in J}) \defeq E(\bigvee_{j \in J} f_j) \defeq \menge{\ve_k}{n_k \neq 0}$. Moreover, given another element $d \in E\reg$, we write $d \cdot \gekl{f_j}_{j \in J} \defeq \menge{d f_j}{j \in J} \eqdef \gekl{f_j}_{j \in J} \cdot d$ and $(d \cdot \gekl{f_j}_{j \in J})\reg \defeq (d \cdot \gekl{f_j}_{j \in J}) \cap E\reg = (\gekl{f_j}_{j \in J} \cdot d) \cap E\reg \eqdef (\gekl{f_j}_{j \in J} \cdot d)\reg$.

Now let $E$ be a semilattice, and let $\Gamma$ be a group acting on $E$ via semigroup automorphisms denoted by $e \ma \tau_g(e)$ ($g \in \Gamma$). Let us assume that we are given a collection of finite covers $\sR$ for $E$, i.e., for every $e \in E\reg$ a set $\sR(e)$ of finite covers for $e$. Let $I_{\Zz}$ be the ideal $\spkl{\menge{e - \bigvee\cR}{e \in E\reg, \ \cR \in \sR(e)}}_{\Zz} \triangleleft \Zz_0[E]$ of $\Zz_0[E]$, and assume that the $\Gamma$-action on $E$ or rather $\Zz_0[E]$ induces a $\Gamma$-action on the quotient $\Zz_0[E] / I_{\Zz}$. Furthermore, let $I$ be the $\Gamma$-invariant ideal $\spkl{\menge{e - \bigvee\cR}{e \in E\reg, \ \cR \in \sR(e)}} \triangleleft C^*_u(E)$ of $C^*_u(E)$. $I$ is the ideal of $C^*_u(E)$ generated by $I_{\Zz}$. Consider the $\Gamma$-action on the quotient $C^*_u(E) / I$ induced by the $\Gamma$-action on $E$.

\btheo
\label{maintheo}

In the situation above, assume that the following conditions are satisfied:
\begin{enumerate}

\item[(i)] For $d$, $e$ in $E\reg$ with $de \neq 0$ and $\cR \in \sR(e)$, either $de \in (d \cdot \cR)\reg$ or $(d \cdot \cR)\reg \in \sR(d e)$.

\item[(ii)] For $e \in E\reg$, pairwise distinct $\cR_1, \dotsc, \cR_r$ in $\sR(e)$ and $\ve_i \in E(\bigvee\cR_i)$ for $1 \leq i \leq r$, we have for every $1 \leq j \leq r$: If 
$\prod_{\substack{i=1 \\ i \neq j}}^r \ve_i \neq 0$, then $\prod_{i=1}^r \ve_i \lneq \prod_{\substack{i=1 \\ i \neq j}}^r \ve_i$. Note that for $r=1$, we set the product $\prod_{\substack{i=1 \\ i \neq j}}^r \ve_i$ as $e$.

\item[(iii)] For every $g \in \Gamma$ and $e \in E\reg$, we have $\tau_g (\sR(e)) = \sR(\tau_g(e))$.

\end{enumerate}
Then Theorem~4.11 in \cite{L-N} gives rise to an algebraic independent resolution of $\Gamma \curvearrowright \Zz_0[E] / I_{\Zz}$, and hence also to an independent resolution of $\Gamma \curvearrowright C^*_u(E) / I$.

If we have, in addition to the assumptions above, that 

\begin{enumerate}
\item[(iv)]$\sup_{e \in E\reg} \abs{\sR(e)} < \infty$,
\end{enumerate}
then the independent resolution of $\Gamma \curvearrowright C^*_u(E) / I$ from above is of length at most $\sup_{e \in E\reg} \abs{\sR(e)}$.
\etheo
\bproof
This is an immediate consequence of \cite[Theorem~4.11]{L-N} and Corollary~\ref{algindres-->Cstarindres}.
\eproof

\section{Quotients of inverse semigroup C*-algebras}
\label{Q-isC}

Reduced C*-algebras of 0-F-inverse semigroups which admit gradings injective on maximal elements (in the sense of \cite{Nor2}) can be described up to Morita equivalence as crossed products of totally disconnected dynamical systems which admit an invariant regular basis. This was observed in \cite{Nor2}. Now we consider quotients of such inverse semigroup C*-algebras, for instance tight C*-algebras of these inverse semigroups. We show that if the quotients are given by relations which satisfy conditions analogous to the ones in Theorem~\ref{maintheo}, then these quotients are Morita equivalent to crossed products which admit finite length independent resolutions. This will be an application of Theorem~\ref{maintheo}.

The general framework for the study of these inverse semigroup C*-algebras and their quotients is given by the notion of partial actions of groups on semilattices. We show that such partial actions can be dilated to ordinary actions on enveloping semilattices. Moreover, relations for our original semilattice satisfying conditions analogous to the ones in Theorem~\ref{maintheo} give rise to relations of the enveloping semilattice which satisfy conditions (i) to (iv) from Theorem~\ref{maintheo} with respect to the dilated action.

Let $E$ be a semilattice, let $E^1$ be $E$ if $E$ already has a unit and the unitalization $E \cup \gekl{1}$ otherwise.
\bdefin
A partial automorphism of $E$ is given by the following data:
\begin{itemize}
\item a projection $d \in E^1$ (the domain)
\item a projection $r \in E^1$ (the range)
\item a semigroup isomorphism $\theta: \ dEd \cong rEr$.
\end{itemize}
\edefin
We will usually write $\theta$ for the partial automorphism.

\bdefin
A partial action $\theta$ of a group $\Gamma$ on $E$ is given by partial automorphisms of $E$,
\bgloz
  \theta_g: \ d(g)Ed(g) \cong r(g)Er(g)
\egloz
one partial automorphism for every group element $g \in \Gamma$, such that we have $d(1) = r(1) = 1$, $\theta_1 = \id_E$ for the identity $1 \in \Gamma$, and $\theta_g \circ \theta_h \leq \theta_{gh}$.
\edefin
This last inequality means the following: By definition, the composition $\theta_g \circ \theta_h$ of $\theta_g$ with $\theta_h$ is given by
\bgloz
  \theta_h^{-1}(r(h)Er(h) \cap d(g)Ed(g)) \to \theta_g(r(h)Er(h) \cap d(g)Ed(g)), \ e \ma \theta_g(\theta_h(e)).
\egloz
Note that $r(h)Er(h) \cap d(g)Ed(g) = (r(h)d(g))E(r(h)d(g))$. Therefore, $\theta_g \circ \theta_h$ is again a partial automorphism of $E$ in our sense, with domain $\theta_h^{-1}(r(h)d(g))$ and range $\theta_g(r(h)d(g))$. We observe that
\bgloz
  \theta_h^{-1}(r(h)d(g))E\theta_h^{-1}(r(h)d(g))
   = \menge{e \in E}{e \leq d(h) \text{ and } \theta_h(e) \leq d(g)}.
\egloz
So the projections in $\theta_h^{-1}(r(h)d(g))E\theta_h^{-1}(r(h)d(g))$ are precisely those projections for which it makes sense to apply $\theta_h$ and then $\theta_g$. The condition $\theta_g \circ \theta_h \leq \theta_{gh}$ means that for every $e \in \theta_h^{-1}(r(h)d(g))E\theta_h^{-1}(r(h)d(g))$, we want to have $\theta_g(\theta_h(e)) = \theta_{gh}(e)$. For this to make sense, we need to have $\theta_h^{-1}(r(h)d(g)) \leq d(gh)$. This is part of the requirement when we ask for the condition $\theta_g \circ \theta_h \leq \theta_{gh}$.

It is obvious that a partial action $\theta$ of $\Gamma$ on $E$ induces in a canonical way a partial action of $\Gamma$ on $C^*_u(E)$, and we again denote this partial action by $\theta$.

Given a partial action $\theta$ of a group $\Gamma$ on $E$, we construct the enveloping semilattice and the dilated action. First, we introduce the following equivalence relation on $\Gamma \times E$:
\bgloz
  (g,d) \sim (h,e) \LRarr \theta_{h^{-1}g}(d)=e.
\egloz
More precisely, the equation $\theta_{h^{-1}g}(d)=e$ includes the requirement that $d \leq d(h^{-1}g)$. It is clear that $\sim$ indeed defines an equivalence relation. The equivalence class of $(h,e)$ will be denoted by $[h,e]$. Moreover, it is easy to check that the formula
\bgloz
  [g,d] \cdot [h,e] \defeq [g,d \theta_{g^{-1}h}(e d(g^{-1}h))]
\egloz
defines a product on $\Gamma \times E / \sim$ so that $(\Gamma \times E / \sim, \cdot)$ becomes a semilattice.

\bdefin
We define a semilattice $\Env(E)$ by setting $\Env(E) \defeq (\Gamma \times E / \sim, \cdot)$.
\edefin

It is easy to see that for every $g \in \Gamma$, the map $[h,e] \ma [gh,e]$ is a well-defined automorphism of $\Env(E)$.
\bdefin
\label{def-tau}
We let $\tau$ be the action of $\Gamma$ on $\Env(E)$ given by $\tau_g [h,e] = [gh,e]$, and we denote the induced $\Gamma$-action on $C^*_u(\Env(E))$ by $\tau$ as well.
\edefin

It is easy to see that the map $E \to \Env(E), \ e \ma [1,e]$ defines an injective homomorphism of semilattices. Moreover, the partial action $\theta$ of $\Gamma$ on $E$ induces a partial action $\theta^1$ of $\Gamma$ on $E^1$ with $\theta^1_g \defeq \theta_g$ if $d(g), r(g) \in E$ and where $\theta^1_g$ is the unique unital extension of $\theta_g$ if $d(g) = r(g) = 1$. Our construction applied to $E^1$ and $\theta^1$ yields another semilattice $\Env(E^1)$ with a $\Gamma$-action. Again $E^1$ sits as a subsemilattice in $\Env(E^1)$. Also, $\Env(E)$ sits canonically as a $\Gamma$-invariant ideal in $\Env(E^1)$. Let $1$ be the unit of $E^1$. Then $1 (\Env(E)) 1$ is the subsemilattice of $\Env(E)$ corresponding to $E$. On the level of C*-algebras, we have that $C^*_u(\Env(E))$ is an essential ideal of $C^*_u(\Env(E^1))$, so that we can think of $1 \in C^*_u(\Env(E^1))$ as a multiplier of $C^*_u(\Env(E))$. In addition, we can canonically identify $C^*_u(E)$ with $1(C^*_u(\Env(E)))1$.

\bremark
It is straightforward to check that $(\tau,C^*_u(\Env(E)))$ is the enveloping action of $(\theta,C^*_u(E))$, in the sense of \cite[Definition~2.3]{Aba}. Therefore, by \cite[Proposition~2.1]{Aba}, the dual action $(\hat{\tau},\widehat{\Env(E)})$ of $(\tau,C^*_u(\Env(E)))$ is the enveloping action of the dual action $(\hat{\theta},\widehat{E})$ of $(\theta,C^*_u(E))$. In particular, $(\hat{\theta},\widehat{E})$ admits an enveloping action on a Hausdorff space.
\eremark

With this remark in mind, the following lemma is not surprising.
\blemma
\label{part-fullcorner}
We have an ismorphism $C^*_u(E) \rtimes_{\theta,r} G \cong 1 \rukl{ C^*_u(\Env(E)) \rtimes_{\tau,r} G } 1$ determined by $e V_g \ma e \cdot (1 U_g 1) = e U_g 1$ for all $e \in E$ and $g \in \Gamma$, and the latter C*-algebra is a full corner of $C^*_u(\Env(E)) \rtimes_{\tau,r} G$. Here $V_g$ is the canonical partial isometry in the multiplier algebra of $C^*_u(E) \rtimes_{\theta,r} G$ corresponding to $g \in \Gamma$, and $U_g$ is the canonical unitary in the multiplier algebra of $C^*_u(\Env(E)) \rtimes_{\tau,r} G$ for $g \in \Gamma$.
\elemma
\bproof
This is an immediate consequence of the observation that $(\tau,C^*_u(\Env(E)))$ is the enveloping action of $(\theta,C^*_u(E))$ and of the construction of reduced crossed products (see for instance \cite{McCl}).

Given a faithful, non-degenerate representation $\pi: \ C^*_u(\Env(E)) \to \cL(H)$, we extend $\pi$ to $C^*_u(\Env(E^1))$ so that we can form $\pi(1)$. Let $\ti{\pi}$ be the twisted representation $C^*_u(\Env(E)) \to \cL(H \otimes \ell^2 \Gamma)$ given by $\ti{\pi}(x)(\xi \otimes \ve_{\gamma}) = \pi(\tau_{\gamma^{-1}}(x))(\xi) \otimes \ve_{\gamma}$. Since $\ti{\pi}$ is again non-degenerate, we can extend it to $C^*_u(\Env(E^1))$ and form $\ti{\pi}(1)$. Let $\lambda$ be the left regular representation of $\Gamma$ on $\ell^2 \Gamma$ and form $1 \otimes \lambda: \ \Gamma \to U(H \otimes \ell^2 \Gamma)$. The reduced crossed product $C^*_u(\Env(E)) \rtimes_{\tau,r} \Gamma$ is by definition the C*-algebra generated by $\ti{\pi}(x) (1 \otimes \lambda_g)$ for $x \in C^*_u(\Env(E))$ and $g \in \Gamma$. Now $\pi \vert_{C^*_u(E)} : \ C^*_u(E) \to \cL(\pi(1) H)$ is a faithful representation of $C^*_u(E)$, and the representation $\rukl{\pi \vert_{C^*_u(E)}}^{\sim}$ (using the notation from \cite[\S~3]{McCl}) is just the cut-down of $\ti{\pi} \vert_{C^*_u(E)}$ by $\ti{\pi}(1)$. Moreover, for every $g \in \Gamma$, $\ti{\pi}(1) (1 \otimes \lambda_g) \ti{\pi}(1)$ is just the partial isometry used in the definition of reduced partial crossed products in \cite[\S~3]{McCl}. The first part of our lemma follows. That $1 \rukl{ C^*_u(\Env(E)) \rtimes_{\tau,r} \Gamma } 1$ is a full corner follows immediately from the obvious fact that $\Env(E) = \bigcup_{g \in \Gamma} \tau_g(E)$.
\eproof

Now let us consider relations.

\blemma
\label{RE-EnvE}
Assume that $\theta$ is a partial action of a group $\Gamma$ on $E$. For every $e \in E\reg$, let $\sR(e)$ be a finite set of finite covers for $e$ such that the following conditions hold:
\begin{itemize}

\item[(1p)] For $d$, $e$ in $E\reg$ with $de \neq 0$ and $\cR \in \sR(e)$, either $de \in (d \cdot \cR)\reg$ or $(d \cdot \cR)\reg \in \sR(d e)$.

\item[(2p)] For $e \in E\reg$, pairwise distinct $\cR_1, \dotsc, \cR_r$ in $\sR(e)$ and $\ve_i \in E(\bigvee\cR_i))$ for $1 \leq i \leq r$, we have for every $1 \leq j \leq r$: If $\prod_{\substack{i=1 \\ i \neq j}}^r \ve_i \neq 0$, then $\prod_{i=1}^r \ve_i \lneq \prod_{\substack{i=1 \\ i \neq j}}^r \ve_i$. As before, we define the product $\prod_{\substack{i=1 \\ i \neq j}}^r \ve_i$ to be $e$ in the case $r=1$.

\item[(3p)] For every $g \in \Gamma$ and $e \in E\reg$ with $e \leq d(g)$, we have $\tau_g (\sR(e)) = \sR(\theta_g(e))$.

\item[(4p)] $\sup_{e \in E\reg} \abs{\sR(e)} < \infty$.

\end{itemize}

If we now set for $[h,e] \in \Env(E)\reg$: $\sR([h,e]) \defeq \tau_h(\sR(e))$, then $\Env(E)$ and $\sR(x)$, $x \in \Env(E)\reg$ satisfy the conditions (i) to (iv) from Theorem~\ref{maintheo}.
\elemma
\bproof
It is easy to see that for every $x \in \Env(E)\reg$, $\sR(x)$ is a well-defined finite set of finite covers for $x$. Moreover, conditions (ii), (iii) and (iv) are easy to check. It remains to check condition (i). Let $x=[g,d]$ and $y=[h,e]$ be elements in $\Env(E)\reg$ with $xy \neq 0$. We have to show that for all $\cR \in \sR(e)$, either $xy$ lies in $x \cdot \tau_h(\cR)$ or $(x \cdot \tau_h(\cR))\reg \in \sR(xy)$. First, let us see that we can without loss of generality assume that $x$ lies in $E$. Namely, $x = \tau_g([1,d])$, and we have $xy = \tau_g([1,d]\tau_g^{-1}y)$, $x \cdot \tau_h(\cR) = \tau_g([1,d] \tau_g^{-1} \tau_h(\cR))$ and $\sR(xy) = \sR(\tau_g([1,d]\tau_g^{-1}y)) = \tau_g(\sR([1,d]\tau_g^{-1}y))$. This means that once we prove our claim for $[1,d]$ in place of $x$ and $\tau_g^{-1}y$ in place of $y$, we are done. In other words, we can assume that $g = 1$.

For $[1,f] \in \cR$, we compute $[1,d]\tau_h[1,f] = [1,d][h,f] = [1,d\theta_h(fd(h))]$. Condition (1p) tells us that either $ed(h) \in \cR \cdot d(h)$ or that $(\cR \cdot d(h))\reg \in \sR(ed(h))$. In the first case, we conclude that $xy = [1,d][h,e] = [1,d\theta_h(ed(h))] = [1,d]\tau_h[1,f]$ (for some $f$) lies in $x \cdot \tau_h(\cR)$. In the second case, it follows that $(\theta_h(\cR \cdot d(h)))\reg \in \sR(\theta_h(ed(h)))$ by condition (3p). Now, condition (1p) again says that we either have $d \theta_h(ed(h)) \in d(\theta_h(\cR \cdot d(h)))$ or $(d(\theta_h(\cR \cdot d(h))))\reg \in \sR(d(\theta_h(ed(h))))$. In the first case, we have $xy = [1,d\theta_h(ed(h))] \in x \cdot \tau_h(\cR)$. In the second case, we conclude that $(x \cdot \tau_h(\cR))\reg \in \sR(xy)$ since $xy = d(\theta_h(ed(h)))$.
\eproof

\bprop
\label{ExG/I}
In the situation of Lemma~\ref{RE-EnvE}, set
\bglnoz
  I &\defeq& \spkl{\menge{e - \bigvee\cR}{e \in E\reg, \ \cR \in \sR(e)}} \ \triangleleft \ C^*_u(E), \\
  \Env(I) &\defeq& \spkl{\menge{x - \bigvee\cR}{x \in \Env(E)\reg, \ \cR \in \sR(x)}} \ \triangleleft \ C^*_u(\Env(E)).
\eglnoz
Let $\spkl{I} \defeq \spkl{I}_{C^*_u(E) \rtimes_{\theta,r} \Gamma}$ be the ideal of $C^*_u(E) \rtimes_{\theta,r} \Gamma$ generated by $I$ and $\spkl{\Env(I)} \defeq \spkl{\Env(I)}_{C^*_u(\Env(E)) \rtimes_{\tau,r} \Gamma}$ be the ideal of $C^*_u(\Env(E)) \rtimes_{\tau,r} \Gamma$ generated by $\Env(I)$.

Then $I = 1 (\Env(I)) 1$, $\spkl{I} = 1 \spkl{\Env(I)} 1$, and the isomorphism from Lemma~\ref{part-fullcorner} induces an isomorphism $\rukl{C^*_u(E) \rtimes_{\theta,r} \Gamma} / \spkl{I} \cong \dot{1} \rukl{ \rukl{C^*_u(\Env(E)) \rtimes_{\tau,r} \Gamma} / \spkl{\Env(I)} } \dot{1}$. If furthermore $\Gamma$ is exact, then the isomorphism from Lemma~\ref{part-fullcorner} also induces an isomorphism $\rukl{C^*_u(E) \rtimes_{\theta,r} \Gamma} / \spkl{I} \cong \dot{1} \rukl{ \rukl{C^*_u(\Env(E)) / \Env(I)} \rtimes_{\tau,r} \Gamma} \dot{1}$. Here $\dot{1}$ is the image of $1$ in the multiplier algebra of the corresponding quotient, and $\dot{1}$ gives rise to a full corner (regardless whether $\Gamma$ is exact or not). In addition, $\Gamma \curvearrowright C^*_u(\Env(E)) / \Env(I)$ admits a finite length independent resolution.
\eprop
\bproof
The equation $I = 1 (\Env(I)) 1$ follows from (1p) and (3p). $\spkl{I} = 1 \spkl{\Env(I)} 1$ is an immediate consequence. The rest follows from Lemma~\ref{part-fullcorner}, Lemma~\ref{RE-EnvE} and Theorem~\ref{maintheo}.
\eproof

Now let $S$ be an inverse semigroup with zero element. For $s \in S$ let $\Lambda(s)$ be the partial isometry on $\ell^2 S\reg$ defined by $\Lambda(s) \ve_x = \ve_{sx}$ if $s^*sx = x$ and $\Lambda(s) \ve_s = 0$ otherwise. By definition, $C^*_r(S)$ is the C*-algebra generated by $\Lambda(s)$, $s \in S$. Note that we consider partial isometries on $\ell^2 S\reg$ to make sure that $\Lambda(0) = 0$.

Let $S$ be 0-F-inverse, and let $G$ be a group and $\sigma: \ (S^1)\reg \to G$ a morphism injective on the set of maximal elements $M(S^1)$ as in \cite[\S~1]{Nor2}. Let $s_g$ be the maximal element of $\sigma^{-1}(g)$ if the latter set is non-empty, and let $s_g \defeq 0$ otherwise. We denote by $E$ the semilattice of idempotent elements in $S$. In such a situation, a partial action $\theta$ of $G$ on $E$ is given as follows: For $g \in G$, we set $d(g) \defeq s_g^* s_g$, $r(g) \defeq s_g s_g^*$ and $\theta_g: d(g)Ed(g) \to r(g)Er(g), \ e \ma s_g e s_g^*$. First of all, let us prove the following
\blemma
\label{ExG-S}
We have an isomorphism $C^*_u(E) \rtimes_{\theta,r} G \cong C^*_r(S)$ determined by $e V_g \ma \Lambda(e s_g)$ for all $e \in E$ and $g \in G$.
\elemma
\bproof
The map $S\reg \to E\reg \times G$, $s \ma (ss^*,\sigma(s))$ is injective since $s = s s^* s_{\sigma(s)}$. Using this map, we view $S\reg$ as a subset of $E\reg \times G$, and we let $P \in \cL(\ell^2 E \reg \otimes \ell^2 G)$ be the orthogonal projection onto $\ell S\reg \subseteq \ell^2 E\reg \otimes \ell^2 G$.

Let $\pi: \ C^*_u(\Env(E)) \to \cL(\ell^2 \Env(E)\reg)$ be the left regular representation of the semilattice $E$ viewed as an inverse semigroup. As before, we extend $\pi$ to $C^*_u(\Env(E^1))$ so that we can form $\pi(1)$. It is clear that we can represent $C^*_u(\Env(E)) \rtimes_{\tau,r} G$ faithfully on $\ell^2 \Env(E)\reg \otimes \ell^2 G$ by sending $x \in \Env(E)$ to $\pi(x) \otimes 1$ and $U_g$ to $T_g \otimes \lambda_g$ for $g \in G$, where $T_g(\ve_x) = \ve_{\tau_g(x)}$ for $x \in \Env(E)\reg$. Using the isomorphism $C^*_u(E) \rtimes_{\theta,r} G \cong 1 \rukl{ C^*_u(\Env(E)) \rtimes_{\tau,r} G } 1$ from Lemma~\ref{part-fullcorner}, we obtain a faithful representation of $C^*_u(E) \rtimes_{\theta,r} G$ on $\pi(1)(\ell^2 \Env(E)\reg) \otimes \ell^2 G = \ell^2 E\reg \otimes \ell^2 G$ given by $e V_g \ma (\pi(e) \otimes 1) (\pi(1) \otimes 1) (T_g \otimes \lambda_g) (\pi(1) \otimes 1)$. An obvious computation shows that both $\pi(e) \otimes 1$, $e \in E$, and $(\pi(1) \otimes 1) (T_g \otimes \lambda_g) (\pi(1) \otimes 1)$, $g \in G$, leave the subspace $\ell^2 S\reg \subseteq \ell^2 E\reg \otimes \ell^2 G$ invariant. Moreover, we have $P (\pi(e) \otimes 1) (\pi(1) \otimes 1) (T_g \otimes \lambda_g) (\pi(1) \otimes 1) P = \Lambda(e s_g)$. Therefore, cutting down by $P$ gives rise to a surjective homomorphism $C^*_u(E) \rtimes_{\theta,r} G \to C^*_r(S)$. This homomorphism is injective since it fits into the following commutative diagram
\bgloz
  \xymatrix@C=1mm{
  C^*_u(E) \rtimes_{\theta,r} G \ar[rr] \ar[d] & \ \ \ \ \ \ & C^*_r(S) \ar[d] \\
  C^*_u(E) \ar[rr]^{\id} & \ \ \ \ \ \ & C^*_u(E)
  }
\egloz
where the vertical arrows are given by the canonical faithful conditional expectations.
\eproof

Combining Lemma~\ref{part-fullcorner} with Lemma~\ref{ExG-S}, we obtain the following
\bcor
\label{S-fullcorner}
$C^*_r(S)$ is isomorphic to the full corner $1 \rukl{ C^*_u(\Env(E)) \rtimes_{\tau,r} G } 1$ of $C^*_u(\Env(E)) \rtimes_{\tau,r} G$ via $\Lambda(e s_g) \ma e \cdot (1 U_g 1) = e U_g 1$.
\ecor
This makes the observations from \cite[\S~2]{Nor2} a bit more explicit.

Again, we turn to relations and the corresponding ideals. The following is an immediate consequence of our discussions:
\bprop
\label{relations-for-S}
Let $\theta$ be the partial action of a group $G$ on a semilattice $E$ attached to a 0-F-inverse semigroup $S$ and a morphism $\sigma$: $(S^1)\reg \to G$ injective on $M(S^1)$ as above. Assume that for every $e \in E\reg$, we are given a finite set $\sR(e)$ of finite covers for $e$ such that conditions (1p) to (4p) from Lemma~\ref{RE-EnvE} hold. Set $I \defeq \spkl{\menge{e - \bigvee\cR}{e \in E\reg, \ \cR \in \sR(e)}} \triangleleft C^*_u(E)$, let $\spkl{I}$ be ideal of $C^*_r(S)$ or $C^*_u(E) \rtimes_{\theta,r} G$, respectively, which is generated by $I$, and let $\Env(I)$, $\spkl{\Env(I)}$ be as in Proposition~\ref{ExG/I}.

If $G$ is exact, then the isomorphisms from Lemma~\ref{ExG-S} and Proposition~\ref{ExG/I} give rise to isomorphisms
$$
  C^*_r(S) / \spkl{I} \cong C^*_u(E) \rtimes_{\theta,r} G / \spkl{I} \cong \dot{1} \rukl{ \rukl{C^*_u(\Env(E)) / \Env(I)} \rtimes_{\tau,r} G} \dot{1}.
$$
The latter C*-algebra is a full corner, so that all these C*-algebras are Morita equivalent to $\rukl{C^*_u(\Env(E)) / \Env(I)} \rtimes_{\tau,r} G$. And finally, $G \curvearrowright C^*_u(\Env(E)) / \Env(I)$ admits a finite length independent resolution.
\eprop

\bremark
The dual system $(\widehat{\Env(E)},G,\hat{\tau})$ in our setting can be canonically identified with the dynamical system $(\Omega,G,\tau)$ from \cite[\S~2]{Nor2}.
\eremark

\bremark
Assume that in Proposition~\ref{relations-for-S}, we can choose relations $\sR(e)$, $e \in E\reg$ in such a way that the spectrum of $C^*_u(E)/I$ identifies with the tight spectrum $\widehat{E}_{\tight}$ in the sense of \cite{Exel}. Then Proposition~\ref{relations-for-S} gives a way to describe the tight (reduced) C*-algebra of $S$ as a crossed product which admits a finite length independent resolution.
\eremark

\section{Computing K-theory in the case of free actions}
\label{k-theory-free-action-section}

Let $E$ be a fixed $\Gamma$-semilattice with a fixed system $\sR$ of covers satisfying (i)-(iii) of Theorem \ref{maintheo}. Suppose also that $\Gamma$ acts freely on $E\reg$. This situation was also discussed in \cite[\S 5]{L-N} where we found methods for computing $H_*(\Gamma,\Zz_0(E)/I_\Zz)$. If we also suppose that $\Gamma$ is exact and satisfies the Baum Connes conjecture with coefficients we can use information about these homology groups to describe the K-theory of $(C^*_u(E)/I)\rtimes_r\Gamma$.

\blemma
\label{k-theory-free-action-lemma}
Continue with the assumptions introduced in the beginning of the section. By applying $K_0$ to the sequence
$$
  \dotso \overset{\phi_3}{\lori} C^*_u(E_2)\rtimes_r\Gamma\overset{\phi_2}{\lori} C^*_u(E_1)\rtimes_r\Gamma \overset{\phi_1}{\lori} C^*_u(E)\rtimes_r\Gamma\to 0
$$
one obtains the chain complex
$$
  C=\left(\dotso \to \Zz_0[\Gamma\setminus E_2] \to \Zz_0[\Gamma\setminus E_1] \to \Zz_0[\Gamma\setminus E] \to 0\right)
$$
defined in \cite[\S 5]{L-N}. Moreover there is a $\Gamma$-equivariant isomorphism $K_0(C^*_u(E)/I)\cong \Zz_0(E)/I_\Zz$, and so we have $H_*(C)\cong H_*(\Gamma,K_0(C^*_u(E)/I))$.
\elemma
\bproof
Let $\phi_k$ denote the $*$-homomorphism $\phi_k:C^*_u(E_k)\rtimes_r\Gamma\to C^*_u(E_{k-1})\rtimes_r\Gamma$ in the long exact sequence. By definition, its restriction to $\Zz_0[E_k]$ is the map induced from the inclusion $E_k\hookrightarrow \Zz_0[E_{k-1}]$. Moreover, \cite[Corollary~3.14]{C-E-L2} gives us that $K_0(C^*_u(E_k)\rtimes_r\Gamma)\simeq\Zz_0[\Gamma\setminus E_k]$. In this isomorphism the $K_0$-class of $e\in E_k$ (identified as an element of $C^*_u(E_k)$) is sent to the class $[e]\in\Zz_0[\Gamma\setminus E_k]$. Thus (omitting the isomorphism) $(\phi_k)_*$ maps $[e]$ to $[f]\in\Zz_0[\Gamma\setminus E_{k-1}]$, where $f\in\Zz_0[E_{k-1}]$ is the inclusion of $e$. Then by definition $[f]=\partial_k([e])$, where $\partial_k:\Zz_0[\Gamma\setminus E_k] \to \Zz_0[\Gamma\setminus E_{k-1}]$ is the k'th boundary map in $C$.

The last statement, that $H_*(C)\cong H_*(\Gamma,K_0(C^*_u(E)/I))$ follows as in \cite[\S 5]{L-N}.
\eproof

\bprop
\label{k-theory-free-action-prop}
Continue with the assumptions introduced in the beginning of the section. Let $n=\sup_{e \in E\reg} \abs{\sR(e)}$. Let $D=C^*_u(E)/I$. Then
\begin{itemize}
  \item If $n=1$, $K_0(D\rtimes_r\Gamma)\cong H_0(\Gamma,K_0(D))$ and $K_1(D\rtimes_r\Gamma)\cong H_1(\Gamma,K_0(D))$.
	\item If $n=2$, $K_0(D\rtimes_r\Gamma)\cong H_0(\Gamma,K_0(D))\oplus H_2(\Gamma,K_0(D))$ and $K_1(D\rtimes_r\Gamma)\cong H_1(\Gamma,K_0(D))$.
	\item If $n=3$ and $H_3(\Gamma,K_0(D))=0$, there is an extension
	$$
	  0\to H_0(\Gamma,K_0(D))\to K_0(D\rtimes_r\Gamma)\to H_2(\Gamma,K_0(D))\to 0,
	$$
	and $K_1(D\rtimes_r\Gamma)\cong H_1(\Gamma,K_0(D))$.
\end{itemize}
\eprop
\bproof
For any map $f:X\to Y$ between sets $X,Y$, let $f^\circ$ denote the restriction $f^\circ:X\to f(X)$. Assume $n>0$ and look at the short exact sequence
$$
  0\to C^*_u(E_n)\rtimes_r\Gamma \overset{\phi_n}{\lori} C^*_u(E_{n-1})\rtimes_r\Gamma \overset{\phi_{n-1}^\circ}{\lori} \mathrm{ker}_{n-2}\to 0
$$
where $\mathrm{ker}_{n-2}=\im\phi_{n-1}=\ker\phi_{n-2}$. Using Lemma~\ref{k-theory-free-action-lemma} we get that this short exact sequence induces the six-term exact sequence
$$
  \xymatrix{
	  \Zz_0[\Gamma\setminus E_n] \ar[r]^{\partial_n} & \Zz_0[\Gamma\setminus E_{n-1}] \ar[r]^{(\phi_{n-1}^\circ)_*} & K_0(\mathrm{ker}_{n-2}) \ar[d] \\
		K_1(\mathrm{ker}_{n-2}) \ar[u]                 & 0 \ar[l]                              & 0 \ar[l]
	}
$$
So $K_0(\mathrm{ker}_{n-2})\cong\coker\partial_n$ with $(\phi_{n-1}^\circ)_*$ being the quotient map, and $K_1(\mathrm{ker}_{n-2})\cong\ker\partial_n=H_n(C)$. Assuming for a moment that $n=1$ we get $\mathrm{ker}_{n-2}=D\rtimes_r\Gamma$. Moreover, $H_0(C)=\coker{\partial_n}$. As shown in Lemma~\ref{k-theory-free-action-lemma}, $H_*(\Gamma,K_0(D))=H_*(C)$, so the first point is proved. Continuing the above computations with $n>1$ we look at the next short exact sequence
$$
  0\to\mathrm{ker}_{n-2} \overset{f_{n-2}}{\lori} C^*_u(E_{n-2})\rtimes_r\Gamma \overset{\phi_{n-2}^\circ}{\lori} \mathrm{ker}_{n-3}\to 0
$$
where $f_{n-2}$ is the inclusion $\ker\phi_{n-2}\hookrightarrow C^*_u(E_{n-2})\rtimes_r\Gamma$. This induces the six-term exact sequence
$$
  \xymatrix{
	  \coker\partial_n \ar[r]^{(f_{n-2})_*} & \Zz_0[\Gamma\setminus E_{n-2}] \ar[r]^{(\phi_{n-2}^\circ)_*} & K_0(\mathrm{ker}_{n-3}) \ar[d] \\
		K_1(\mathrm{ker}_{n-3}) \ar[u]       & 0 \ar[l]                     & H_n(C) \ar[l]
	}
$$
Since $\phi_{n-1}=f_{n-2}\phi_{n-1}^\circ$ we get $\partial_{n-1}(x)=(f_{n-2})_*(\phi_{n-1}^\circ)_*(x)=(f_{n-2})_*(x+\im\partial_n)$. This gives us $K_1(\mathrm{ker}_{n-3})\cong\ker (f_{n-2})_*=\ker\partial_{n-1}/\im\partial_n=H_{n-1}(C)$. As $H_n(C)=\ker\partial_n$ is free over $\Zz$ we get $K_0(\mathrm{ker}_{n-3})\cong H_n(C)\oplus\coker (f_{n-2})_*=H_n(C)\oplus\coker\partial_{n-1}$. Here $(\phi_{n-2}^\circ)_*$ is the quotient map onto $\coker\partial_{n-1}$. If we assume for a moment that $n=2$, then $\mathrm{ker}_{n-3}=D\rtimes_r\Gamma$. Moreover, $\coker\partial_{n-1}=H_0(C)$, so the second point is proved. Continuing the computations for $n>2$ we get the short exact sequence
$$
  0\to\mathrm{ker}_{n-3} \overset{f_{n-3}}{\lori} C^*_u(E_{n-3})\rtimes_r\Gamma \to \mathrm{ker}_{n-4}\to 0
$$
and the associated six-term exact sequence
$$
  \xymatrix{
	  H_n(C)\oplus\coker\partial_{n-1} \ar[r]^{(f_{n-3})_*} & \Zz_0[\Gamma\setminus E_{n-3}] \ar[r] & K_0(\mathrm{ker}_{n-4}) \ar[d] \\
		K_1(\mathrm{ker}_{n-4}) \ar[u]                      & 0 \ar[l]                              & H_{n-1}(C) \ar[l]
	}
$$
Using a similar argument as for $(f_{n-2})_*$ we get that $(f_{n-3})_*(x,y+\im\partial_{n-1})=g(x)+\partial_{n-2}(y)$ for some map $g$. Now if $H_n(C)=0$, $K_1(\mathrm{ker}_{n-4})\cong\ker (f_{n-3})_*= H_{n-2}(C)$. We also see that there is an extension
$$
0\to \coker (f_{n-3})_*\to K_0(\mathrm{ker}_{n-4})\to H_{n-1}(C)\to 0 
$$
If $H_n(C)=0$, then $\coker (f_{n-3})_*=\coker\partial_{n-2}$. In particular, if $n=3$ and $H_3(C)=0$, $\coker (f_{n-3})_*=H_0(C)$. This finishes the proof.
\eproof

\bremark\label{chain-complex-remark}
As noted in Lemma~\ref{k-theory-free-action-lemma}, the homology groups $H_*(\Gamma,K_0(C^*_u(E)/I))$ may be computed as the homology groups of the chain complex $C$ of \cite[\S~5]{L-N}. If the system of covers $\sR$ also satisfies the conditions (A)-(C) of \cite[\S 5]{L-N}, one may due to \cite[Remark~5.6]{L-N} replace the chain complex $C$ with the chain complex $\wt{C}$ in the situation of that remark. The chain complex $\wt{C}$ is also defined in \cite[\S~5]{L-N}.
\eremark

\section{Examples}

\subsection{Graph C*-algebras}

We show that using independent resolutions, it is easy to compute K-theory for graph C*-algebras. We use the same notation as in \cite[\S~5]{Nor2}: Let $\cE = (\cE^0, \cE^1, \sigma, \rho)$ be a graph, and let $S_{\cE}$ be its graph inverse semigroup. $S_{\cE}^1$ is strongly 0-F-inverse with universal grading $(S_{\cE}^1)\reg \to \Fz$ , where $\Fz$ is the free group on $\cE^1$. The semilattice $E$ of idempotent elements in $S_{\cE}$ can be identified with $\cE^* \cup \gekl{0}$, where $\cE^*$ is the set of finite paths of $\cE$. Multiplicaton in $E$ is given by $\mu \cdot \nu \defeq \mu$ if $\nu = \mu \nu'$ for some $\nu' \in \cE^*$, $\mu \cdot \nu \defeq \nu$ if $\mu = \nu \mu'$ for some $\mu' \in \cE^*$, and $\mu \cdot \nu \defeq 0$ otherwise. Here $\mu \nu'$ stands for concatenation of paths.

The partial action of $\Fz$ on $E$ attached to $S_{\cE}$ in the sense of \S~\ref{Q-isC} is given as follows: We view $\cE^*$ as a subset of $\Fz$ in a canonical way. For paths $\mu$ and $\nu$ in $\cE^*$ with length at least one and $\sigma(\mu) = \sigma(\nu)$, let $d(\mu \nu^{-1}) = \nu \cdot E$, $r(\mu \nu^{-1}) = \mu \cdot E$ and $\theta_{\mu \nu^{-1}}(\nu \cdot \xi) \defeq \mu \cdot \xi$. For the identity $1 \in \Fz$, we set $\theta_1 \defeq \id_E$, and all the remaining $g \in \Fz$ do not lie in the image of our grading.

As observed in \cite[\S~5]{Nor2}, $C^*_r(S_{\cE})$ is canonically isomorphic to the Toeplitz algebra of $\cE$. The graph C*-algebra of $\cE$ is the tight version of $C^*_r(S_{\cE})$, i.e., a quotient by a certain ideal. To describe the graph C*-algebra of $\cE$, we consider the following relations: Let $\cE^0_0$ be the vertices $v$ of $\cE$ for which $0 < \# \menge{\kappa \in \cE^1}{v = \rho(\kappa)} < \infty$, and let $\sigma^{-1}(\cE^0_0)$ be the set of paths $\mu$ with $\sigma(\mu) \in \cE^0_0$. For every $\mu \in \sigma^{-1}(\cE^0_0)$, we let $\cR(\mu)$ be the finite cover $\menge{\mu \kappa}{\kappa \in \cE^1, \, \sigma(\mu) = \rho(\kappa)}$ for $\mu$, and we set $\sR(\mu) \defeq \gekl{\cR(\mu)}$. For the remaining $\mu \in E\reg$ which are not in $\sigma^{-1}(\cE^0_0)$, just set $\sR(\mu) \defeq \emptyset$. Let $e_{\mu}$ be the projection in $C^*_u(E)$ corresponding to $\mu \in E$. If we now set
$$
I \defeq \spkl{\menge{e_{\mu} - \bigvee_{\nu\in\cR(\mu)}e_{\nu}}{\mu \in \sigma^{-1}(\cE^0_0)}} \triangleleft C^*_u(E),
$$
then it is clear by construction that the graph C*-algebra $C^*(\cE)$ is canonically isomorphic to $C^*_r(S_{\cE})/\spkl{I}$. Moreover, it is easy to see that the partial action $\Fz \curvearrowright E$ and $\sR(\mu)$, $\mu \in E\reg$, satisfy conditions (1p) to (4p) from Lemma~\ref{RE-EnvE}. 

Proposition~\ref{S-fullcorner} implies that $C^*_r(S_{\cE})$ (and hence the Toeplitz algebra of $\cE$) is isomorphic to a full corner in $C^*_u(\Env(E)) \rtimes_r \Fz$, and Proposition~\ref{relations-for-S} implies that $C^*_r(S_{\cE})/\spkl{I}$ (and hence the graph C*-algebra of $\cE$) is isomorphic to a full corner in $(C^*_u(\Env(E)) / \Env(I)) \rtimes_r \Fz$.

Let us now come to K-theory. Since the stabilizer groups for the action of $\Fz$ on $\Env(E)\reg$ are trivial (see \cite{Nor2}) we could utilize Proposition~\ref{k-theory-free-action-prop}, but in this case it is more illuminating to do the computations directly to illustrate what goes on. Lemma~\ref{RE-EnvE} and \cite[Proposition~4.1]{L-N} yield the semilattice
\bgloz
  E_1 \defeq \menge{[g,\mu] - \bigvee\cR([g,\mu])}{g \in \Fz, \, \mu \in \sigma^{-1}(\cE^0_0)} \cup \gekl{0} \subseteq \Proj(C^*_u(\Env(E))),
\egloz
where $\cR([g,\mu]) = \tau_g(\cR(\mu))$ ($\tau$ is defined in Definition~\ref{def-tau}). By Theorem~\ref{maintheo}, since we have $\sup_{e \in E\reg} \abs{\cR(e)} = 1$, we obtain a short exact sequence
\bgl
\label{graph-ses}
  0 \to C^*_u(E_1) \rtimes_r \Fz \overset{i}{\lori} C^*_u(\Env(E)) \rtimes_r \Fz \to (C^*_u(\Env(E)) / \Env(I)) \rtimes \Fz \to 0.
\egl
With \cite[Corollary~3.14]{C-E-L2} (see also \cite{Nor2}), we compute $K_*(C^*_u(\Env(E)) \rtimes_r \Fz) \cong \oplus_{v \in \cE^0} K_*(\Cz)$, with generators for $K_0$ given by $[e_v]$, $v \in \cE^0$, and $K_*(C^*_u(E_1) \rtimes_r \Fz) \cong \oplus_{w \in \cE^0_0} K_*(\Cz)$, with generators for $K_0$ given by $[e_w - \bigvee_{\nu \in \cR_w} e_{\nu}]$, $w \in \cE^0_0$. Using $\bigvee\cR([g,\mu]) = \sum_{\substack{\kappa \in \cE^1 \\ \sigma(\mu) = \rho(\kappa)}} [g,\mu \kappa]$, we obtain that $i_*: \ K_0(C^*_u(E_1) \rtimes_r \Fz) \to K_0(C^*_u(\Env(E)) \rtimes_r \Fz)$ sends $[e_w - \bigvee_{\nu \in \cR_w} e_{\nu}]$ to $[e_w] - \sum_{\substack{\kappa \in \cE^1 \\ w = \rho(\kappa)}} [e_{\sigma(\kappa)}]$.
Thus we see that $i_*$ can be described using the vertex matrix $A_{\cE}$, i.e., the $\cE^0 \times \cE^0$ matrix given by $A_{\cE}(v,w) = \# \menge{\kappa \in \cE^1}{\rho(\kappa) = v, \, \sigma(\kappa) = w} \in \Nz_0 \cup \gekl{\infty}$. Under the decomposition $\cE^0 = \cE^0_0 \cup (\cE^0 \setminus \cE^0_0)$, $A_{\cE}$ is of the form
$
  \rukl{
  \begin{smallmatrix}
  A_0 & A_1 \\
  * & *
  \end{smallmatrix}
  }
$
where the entries in $*$ are $0$ or $\infty$. Using $A_0$ and $A_1$ from the vertex matrix, $i_*$ identifies with the homomorphism $\begin{bmatrix}I - A_0^t\\ - A_1^t\end{bmatrix}: \ \Zz^{\cE^0_0} \to \Zz^{\cE^0}$. Plugging this result into the six-term exact sequence attached to \eqref{graph-ses}, we obtain for the K-theory of the graph C*-algebra $C^*(\cE)$:
\bgloz
  K_0(C^*(\cE)) \cong \coker\begin{bmatrix}I - A_0^t\\ - A_1^t\end{bmatrix} \text{ and } K_1(C^*(\cE)) \cong \ker\begin{bmatrix}I - A_0^t\\ - A_1^t\end{bmatrix}.
\egloz
This reproves \cite[Theorem~3.1]{D-T}. Note that the chain complex
$$
0 \to \Zz^{\cE^0_0} \overset{i_*}{\lori} \Zz^{\cE^0} \to 0
$$
is easily identified with the chain complex $C$ discussed in \S~\ref{k-theory-free-action-section}, with $i_*=\partial_1$.

\bquestion
Is a similar analysis possible for higher rank graph C*-algebras?
\equestion

\subsection{C*-algebras of one dimensional tilings}

We will see how independent resolutions can be used to compute the K-theory of the C*-algebras associated to one dimensional tilings. A tile in $\Rz$ is a closed interval. A tiling $T$ of $\Rz$ is a set of tiles with pairwise disjoint interiors and union $\Rz$. As in \cite[\S~4.2]{Ke-Law} we describe the connected tiling inverse semigroup as the inverse semigroup associated to a factorial language. Let $\Sigma$ be a finite alphabet (i.e. a finite set). A language $L$ on $\Sigma$ is factorial if for every $x\in L$ every substring of $x$ also belongs to $L$. Assume also that every element of $\Sigma$ occurs in $L$. In our setting we imagine $T$ as a bi-infinite string on a finite set $\Sigma$ of prototiles and $L$ as the factorial language consisting of all finite substrings of $T$.

Let $S_L$ be the inverse semigroup associated to the factorial language $L$. Then the semilattice $E$ of idempotent elements in $S_L$ consists of $0$ as well as all strings on the alphabet $\Sigma\cup\{\check{a}:a\in \Sigma\}$ on the form $x\check{a}y$ where $x,y\in \Sigma^*$ and $xay\in L$. In other words, the nonzero elements of $E$ are elements of $L$ with a check above one of its letters. Multiplication is defined as follows: Let $e,d\in E\reg$ and place $e$ above $d$ such that the checked letter of $e$ is above the checked letter of $d$. If they match on the overlap, glue $e$ and $d$ together on their overlap. If the resulting element belongs to $E$ define this element to be $e\cdot d$. Otherwise $e\cdot d$ is defined to be $0$.

It was shown in \cite{Ke-Law} that $S_\Sigma$ is strongly 0-F-inverse with universal grading $(S_\Sigma^1)^\times\to\Fz$ where $\Fz$ is the free group on the set $\{(a,b)\in\Sigma\times\Sigma: ab\in L\}$. For higher dimensional tilings the connected tiling semigroup is in general not 0-F-inverse. The partial action $\Fz\curvearrowright E$ in the sense of \S~\ref{Q-isC} becomes as follows: With $g\in\Fz$ on the form $g=(a_1,a_2)(a_2,a_3)\cdots (a_{n-2},a_{n-1})(a_{n-1},a_n)$, $a_1,\ldots,a_n\in\Sigma$ we get 
\bglnoz
&& d(g)=\menge{xa_1 a_2\cdots a_{n-1}\check{a_n} y}{x,y\in\Sigma^*,xa_1 a_2\cdots a_{n-1}a_n y\in L}\\
&& r(g)=\menge{x\check{a_1} a_2\cdots a_{n-1}a_n y}{x,y\in\Sigma^*,xa_1 a_2\cdots a_{n-1}a_n y\in L}\\
&& \theta_g(xa_1 a_2\cdots a_{n-1}\check{a_n} y)=x\check{a_1} a_2\cdots a_{n-1}a_n y
\eglnoz
Moreover, $\theta_{g^{-1}}=\theta_g^{-1}$ and $\theta_1=\id_E$. No other $g$ lies in the image of the grading.

For each $e\in E\reg$ set $\cR_1(e)\defeq\menge{ae}{a\in\Sigma, ae\in E}$, $\cR_2(e)\defeq\menge{ea}{a\in\Sigma, ea\in E}$ and let $\sR(e)\defeq\{\cR_1(e),\cR_2(e)\}$. These covers are chosen to make $C^*_r(S_T)/\langle I\rangle$ isomorphic to the tiling C*-algebra $A_T$ of \cite{Ke-Pu}. It is easy to see that the partial action $\Fz \curvearrowright E$ and $\sR(e)$, $e \in E\reg$, satisfy conditions (1p) to (4p) from Lemma~\ref{RE-EnvE}. Since every $x\in L$ is a substring of $T$ we have that $ax\in L$ and $xb\in L$ for at least one $a\in\Sigma$ and one $b\in\Sigma$. Thus $|\sR(e)|=2$ for each $e\in E\reg$.

Let $p(e)\in C^*_u(\Env(E))$ stand for the projection corresponding to $e\in E$ and similarly let $p([g,e])$ stand for the projection corresponding to $[g,e]$ where $g\in\Fz$.

We get the semilattice $E_1$ consisting of $0$ and the elements
\bglnoz
&& p([g,e]||1)\defeq p([g,e])-\bigvee\cR_1([g,e]) = p([g,e])-\sum_{a\in\Sigma,ae\in E}p([g,ae]) \\
&& p([g,e]||2)\defeq p([g,e])-\bigvee\cR_2([g,e]) = p([g,e])-\sum_{a\in\Sigma,ea\in E}p([g,ea])\\
&& p([g,e]||1,2)\defeq p([g,e]||1)p([g,e]||2)
\eglnoz
for each $g\in\Fz$ and $e\in E$. We also get the semilattice $E_2$ consisting of $0$ and the elements
\bglnoz
&& p([g,e]||1|2)\defeq p([g,e]||1) - p([g,e]||1,2) - \sum_{a\in\Sigma, ea\in E} p([g,ea]||1)\\
&& p([g,e]||2|1)\defeq p([g,e]||2) - p([g,e]||1,2) - \sum_{a\in\Sigma, ea\in E} p([g,ae]||2)\\
\eglnoz
for each $g\in\Fz$ and $e\in E$. 

Theorem~\ref{maintheo} gives an exact sequence
\bglnoz
  0 &\to& C^*_u(E_2) \rtimes_r \Fz \to C^*_u(E_1) \rtimes_r \Fz \\
  &\to& C^*_u(\Env(E)) \rtimes_r \Fz \to (C^*_u(\Env(E)) / \Env(I)) \rtimes_r \Fz \to 0.
\eglnoz
As seen in \cite{Nor2}, $\Fz$ acts freely on $\Env(E)\reg$, so with \cite[Corollary~3.14]{C-E-L2} we compute $K_*(C^*_u(\Env(E)) \rtimes_r \Fz) \cong \oplus_{x\in L}K_*(\Cz)$, with generators for $K_0$ given by $[e]$, where $e\in E$ and the first letter of $e$ is checked. Similarly we compute $K_*(C^*_u(E_1) \rtimes_r \Fz)\cong\oplus_{x\in L}(K_*(\Cz)\oplus K_*(\Cz)\oplus K_*(\Cz))$ with generators for $K_0$ given by $[p(e||1)],[p(e||2)],[p(e||1,2)]$ and $K_*(C^*_u(E_2) \rtimes_r \Fz)\cong\oplus_{x\in L}(K_*(\Cz)\oplus K_*(\Cz))$ with generators for $K_0$ given by $[p(e||1|2)],[p(e||2|1)]$ where $e\in E$ and the first letter of $e$ is checked. Applying $K_0$ to our long exact sequence we thus get the chain complex
$$
C = \left(0 \to \bigoplus_{x\in L}(\Zz\oplus\Zz) \overset{\partial_2}{\lori}  \bigoplus_{x\in L}(\Zz\oplus\Zz\oplus\Zz) \overset{\partial_1}{\lori}  \bigoplus_{x\in L} \Zz \to 0\right)
$$
We can now apply Proposition~\ref{k-theory-free-action-prop}. We get $\ker \partial_2=0$, and $H_1(C)\cong\Zz$, generated by $(\sum_{a\in\Sigma}([p(\check{a}||1)]-[p(\check{a}||2)])) + \im \partial_1$. Let $1_x$ be the generator of the $x$'th copy of $\Zz$ in $\bigoplus_L \Zz$ and let $H$ be the subgroup of $\bigoplus_L \Zz$ generated by
$$\menge{1_x-\sum_{a\in\Sigma,ax\in L}1_{ax}}{x\in L}\cup\menge{1_x-\sum_{a\in\Sigma,xa\in L}1_{xa}}{x\in L}.$$
We then get $K_0(A_T)\cong\coker\partial_1\cong (\bigoplus_L \Zz) / H$ and $K_1(A_T)\cong H_1(C)\cong\Zz$.

With some work one can see that this is an affirmation of the observations about the K-theory of one-dimensional tiling C*-algebras found in \cite{Ke-Law}.

\subsection{Boundary quotients of semigroup C*-algebras}

\subsubsection{Right-angled Artin monoids}
\label{A-T}

Let $P$ be a right-angled Artin monoid and $G$ the corresponding Artin group. We refer to \cite{Cr-La1} and \cite{Cr-La2} for details. It is known that $P$ embeds as a subsemigroup into $G$. The left inverse hull $I_l(P)$ is an inverse semigroup of the type studied in \S~\ref{Q-isC}. The corresponding semilattice is given by $\cJ = \menge{pP}{p \in P} \cup \gekl{\emptyset}$ with intersection as multiplication, and the partial action $\theta$ of $G$ on $\cJ$ attached to $I_l(P)$ in \S~\ref{Q-isC} is given by $d(g) = (g^{-1} \cdot P) \cap P \in \cJ$, $r(g) = P \cap (g \cdot P) \in \cJ$ and $\theta_g: \ d(g) \cJ d(g) \to r(g) \cJ r(g)$, $X \ma g \cdot X$. We have canonical isomorphisms $C^*_r(P) \cong C^*_r(I_l(P)) \cong C^*_u(\cJ) \rtimes_{\theta,r} G$. Here $C^*_r(P)$ is the semigroup C*-algebra of $P$, discussed in \cite{{Li1},{Li2}} in a general context and in \cite{{Cr-La1},{Cr-La2}} in the particular case of Artin monoids.

Let us now assume that the underlying graph of our right-angled Artin monoid $P$ is irreducible and finite. Let $S$ be the set of generators of $P$ corresponding to the edges of the graph. In this situation, let us describe the boundary quotient of $C^*_r(P)$ with the help of relations. For each $p \in P$, let $\cR(pP)$ be the finite cover $\menge{psP}{s \in S}$ for $pP \in \cJ$, and set $\sR(pP) \defeq \gekl{\cR(pP)}$. With $I \defeq \spkl{\menge{e_{X} - \bigvee_{Y \in \cR(X)} e_Y}{X \in \cJ\reg}} \triangleleft C^*_u(\cJ)$, \cite[Lemma~3.8 and Corollary~6.6]{Cr-La2} tell us that $C^*_r(P) / \spkl{I}$ is the boundary quotient of $C^*_r(P)$. Moreover, the partial action $\theta$ of $G$ on $\cJ$ and the relations $\sR(X)$, $X \in \cJ\reg$, satisfy conditions (1p) to (4p) from Lemma~\ref{RE-EnvE}: Conditions (2p), (3p) and (4p) are obviously satisfied. Condition (1p) also holds because given $p$, $q$ and $x$ in $P$ with $pP \cap qP = xP$, we have that $x \in pP = \gekl{p} \cup \bigcup_{s \in S} psP$. If $x$ lies in $psP$ for some $s \in S$, then $psP \cap qP = psP \cap pP \cap qP = xP$, and if $x = p$, then $pP \subseteq qP$, thus $psP \cap qP = psP$ for all $s \in S$. The enveloping semilattice of $\cJ$ is given by $\cJ_{P \subseteq G} = \menge{gP}{g \in G} \cup \gekl{\emptyset}$, and $G$ acts by left multiplication. Setting $\cR(gP) \defeq \menge{gsP}{s \in S}$ and $\sR(gP) \defeq \gekl{\cR(gP)}$, Lemma~\ref{RE-EnvE} tells us that $G \curvearrowright \cJ_{P \subseteq G}$ and $\sR(Y)$, $Y \in \cJ_{P \subseteq G}\reg$, satisfy conditions (i) to (iv) of Theorem~\ref{maintheo}. Since $G$ is exact by \cite{G-N}, Proposition~\ref{S-fullcorner} and Proposition~\ref{relations-for-S} imply that $C^*_r(P) \sim_M C^*_u(\cJ_{P \subseteq G}) \rtimes_r G$ and $C^*_r(P) / \spkl{I} \sim_M (C^*_u(\cJ_{P \subseteq G}) / \Env(I)) \rtimes_r G$. \cite[Theorem~4.11]{L-N} yields the semilattice $E_1 = \menge{e_{gP} - \bigvee_{Y \in \cR(gP)} e_Y}{g \in G} \cup \gekl{0} \subseteq \Proj(C^*_u(\cJ_{P \subseteq G}))$, and since $\sup_{Y \in \cJ_{P \subseteq G}\reg} \abs{\sR(Y)} = 1$, we obtain the short exact sequence
\bgloz
  0 \to C^*_u(E_1) \rtimes_r G \overset{i}{\lori} C^*_u(\cJ_{P \subseteq G}) \rtimes_r G \to (C^*_u(\cJ_{P \subseteq G}) / \Env(I)) \rtimes_r G \to 0.
\egloz
\cite[Corollary~3.14]{C-E-L2} yields $K_*(C^*_u(E_1) \rtimes_r G) \cong K_*(\Cz)$, with the generator of $K_0$ given by $[e_P - \bigvee_{s \in S} e_{sP}]$, and $K_*(C^*_u(\cJ_{P \subseteq G}) \rtimes_r G) \cong K_*(\Cz)$, where the generator of $K_0$ is given by $[e_P]$. Since $i_*$ sends $[e_P - \bigvee_{s \in S} e_{sP}]$ to $\chi \cdot [e_P]$, where $\chi$ is the Euler characteristic of the underlying graph of $P$ in the sense of \cite{Cr-La2} and \cite{Iva}, we obtain for the K-theory of the boundary quotient $C^*_r(P) / \spkl{I}$
\begin{itemize}
\item if $\chi = 0$: $K_0(C^*_r(P) / \spkl{I}) \cong K_0((C^*_u(\cJ_{P \subseteq G})/\Env(I)) \rtimes_r G) \cong \Zz$ and \\
$K_1(C^*_r(P) / \spkl{I}) \cong K_1((C^*_u(\cJ_{P \subseteq G})/\Env(I)) \rtimes_r G) \cong \Zz$,
\item if $\chi \neq 0$: $K_0(C^*_r(P) / \spkl{I}) \cong K_0((C^*_u(\cJ_{P \subseteq G})/\Env(I)) \rtimes_r G) \cong \Zz / \abs{\chi}\Zz$ and $K_1(C^*_r(P) / \spkl{I}) \cong K_1((C^*_u(\cJ_{P \subseteq G})/\Env(I)) \rtimes_r G) \cong \gekl{0}$.
\end{itemize} 
We point out that the K-theory of the boundary quotient has already been computed in \cite{Iva} using different methods.

\subsubsection{Group C*-algebras as boundary quotients of semigroup C*-algebras}
\label{groupC}

Under the same assumptions as in \cite[\S~6]{L-N}, we obtain independent resolutions for group C*-algebras. In special cases, for instance in the situation of \cite[\S~6.2]{L-N}, these resolutions have finite length and hence can be used to compute K-theory for group C*-algebras of particular groups.

\subsubsection{Ring C*-algebras for rings of integers}
\label{ringC}

We consider the same partial action $\theta: \ G \curvearrowright \cJ$ as in \S~\ref{min-prim-id}. Let $\cP$ be the set of non-zero prime ideals of $R$. Consider the following relations:
$$\sR((r+\mfa) \times \mfa\reg) = \menge{\menge{(r + s + \mfp \cdot \mfa) \times \mfa\reg}{s \in \mfa / \mfp \cdot \mfa}}{\mfp \in \cP}.$$
With $I \defeq \spkl{\menge{e_{X} - \bigvee_{Y\in\cR}e_Y}{X \in \cJ\reg, \, \cR \in \sR(X)}} \triangleleft C^*_u(\cJ)$, $C^*(R \rtimes R\reg) / \spkl{I}$ is the boundary quotient of $C^*(R \rtimes R\reg)$, hence isomorphic to the ring C*-algebra of $R$ from \cite{Cu-Li1}. It is straightforward to see that $\theta: \ G \curvearrowright \cJ$ and $\sR((r+\mfa) \times \mfa\reg)$, $(r+\mfa) \times \mfa\reg \in \cJ\reg$, satisfy conditions (1p) to (3p) from Lemma~\ref{RE-EnvE}, but (4p) does not hold because $\cP$ is infinite. This problem can be solved as follows: Enumerate the prime ideals, i.e., write $\cP = \gekl{\mfp_1, \mfp_2, \mfp_3, \dotsc}$ and set $\cP_n \defeq \gekl{\mfp_1, \dotsc, \mfp_n}$. Moreover, set $\sR^{(\cP_n)}((r+\mfa) \times \mfa\reg) \defeq \menge{\menge{(r + s + \mfp \cdot \mfa) \times \mfa\reg}{s \in \mfa / \mfp \cdot \mfa}}{\mfp \in \cP_n}$. In this way, we have enforced the finiteness condition (4p), and all the remaining conditions are still satisfied. Let $I^{(\cP_n)}$ be the ideal $\spkl{\menge{e_{X} - \bigvee_{Y\in\cR}e_Y)}{X \in \cJ\reg, \, \cR \in \sR^{(\cP_n)}(X)}}$ of $C^*_u(\cJ)$ corresponding to $\cP_n$. The quotient $C^*(R \rtimes R\reg) / \spkl{I}$ can be identified with the inductive limit $\ilim_n C^*(R \rtimes R\reg) / \spkl{I^{(\cP_n)}}$. Therefore, by continuity of K-theory, it suffices to understand the K-theory of $C^*(R \rtimes R\reg) / \spkl{I^{(\cP_n)}}$. Again, we may apply our results in \S~\ref{alg-Cstalg} and \S~\ref{Q-isC} and proceed as in the previous examples. Although this in principle leads to the K-theory of ring C*-algebras, there are lots of extension problems to be solved along the way, which makes this approach very complicated. Recently, the K-theory for such ring C*-algebras has been completely determined in \cite{Cu-Li2}, \cite{Cu-Li3} and \cite{L-L}, but these computations follow a different route. The key role is played by the so-called duality theorem from \cite{Cu-Li2}.

In a similar fashion, one can also treat the Bost-Connes algebra from \cite{B-C}. However, as far as we can see, this approach does not give a direct computation of the K-theory of the Bost-Connes algebra, unless there is a good understanding of the group homology $H_n(\Qz\pos, \; K_0(C_0(\Az_f))) \cong H_n(\Qz\pos, \; C_0(\Az_f, \Zz))$.

\subsection{Minimal non-zero primitive ideals of $C^*(R \rtimes R\reg)$ and their quotients}
\label{min-prim-id}

Let $K$ be a number field with ring of integers $R$. Consider the $ax+b$-semigroup $P = R \rtimes R\reg$, which is a subsemigroup of $G = K \rtimes K\reg$. Again, consider the partial action $\theta: \ G \curvearrowright \cJ$ attached to the left inverse hull of $P$ as in \S~\ref{A-T}. We have $\cJ = \menge{(r+\mfa) \times \mfa\reg}{r \in R, \; (0) \neq \mfa \triangleleft R} \cup \gekl{\emptyset}$. We view $\cJ$ as a semilattice with multiplication given by intersection of sets. For a non-zero prime ideal $(0) \neq \mfp$ of $R$, let $\cR((r+\mfa) \times \mfa\reg)$ be the finite cover $\menge{(r + s + \mfp \cdot \mfa) \times \mfa\reg}{s \in \mfa / \mfp \cdot \mfa}$ for $(r+\mfa) \times \mfa\reg \in \cJ\reg$, and set $\sR((r+\mfa) \times \mfa\reg) \defeq \gekl{\cR((r+\mfa) \times \mfa\reg)}$. The ideal $I_{\mfp} \defeq \spkl{\menge{e_{(r+\mfa) \times \mfa\reg} - \bigvee_{Y \in \cR((r+\mfa) \times \mfa\reg)} e_Y}{(r+\mfa) \times \mfa\reg \in \cJ\reg}}$ is the minimal non-zero primitive ideal of $C^*_r(P)$ attached to $\mfp$. $\theta: \ G \curvearrowright \cJ$ and $\sR((r+\mfa) \times \mfa\reg)$, $(r+\mfa) \times \mfa\reg \in \cJ\reg$, satisfy conditions (1p) to (4p) from Lemma~\ref{RE-EnvE}. This is proven in \cite[Lemma~3.5]{Li3}, but in a slightly different language. Using our results in \S~\ref{alg-Cstalg} and \S~\ref{Q-isC}, the same procedure as in \S~\ref{A-T} gives a description of the quotient $C^*_r(P) / I_{\mfp}$ as a full corner in a (reduced) crossed product which admits an independent resolution of length one. The corresponding six-term exact sequence can be used to study K-theory. This is worked out in detail in \cite{Li3}, where these ideas lead to a classification result for the semigroup C*-algebras $C^*_r(R \rtimes R\reg)$.

\subsection{The multiplicative boundary quotient of the C*-algebra of $\Nz \rtimes Q$}

A similar, but easier example as in \S~\ref{ringC} is the following: Let $p_1$, $p_2$, ... be the prime numbers (in any order). For a given $n\geq 1$, set $Q = \leck p_1, \dotsc, p_n \rsp$ to be the multiplicative semigroup generated by $p_1,\ldots, p_n$. We form the semidirect product $P \defeq \Nz \rtimes Q$ with respect to the multiplicative action of $Q$ on $\Nz = \gekl{0, 1, 2, \dotsc}$. We set $G \defeq \Zz[\tfrac{1}{p_1}, \dotsc, \tfrac{1}{p_n}] \rtimes \spkl{p_1, \dotsc, p_n}$ and consider the partial action $\theta: \ G \curvearrowright \cJ$ as in \S~\ref{A-T}. $\cJ$ is given by $\menge{(j+m\Nz) \times mQ}{j \in \Nz, \, m \in Q}\cup\{\emptyset\}$. We introduce the relations 
$$
\cR_i((j+m\Nz)\times mQ)\defeq \menge{(j + m r + m p_i \Nz) \times m p_i Q}{0 \leq r \leq p_i-1}
$$
and $\sR((j+m\Nz) \times mQ) \defeq \{\cR_i((j+m\Nz)\times mQ)\}_{i=1}^n$. Let
$$
I \defeq \spkl{\menge{e_X - \bigvee_{Y\in\cR}e_Y}{X \in \cJ\reg, \, \cR \in \sR(X)}} \, \triangleleft \, C^*_u(\cJ)
$$
be the corresponding ideal. A similar analysis as in the previous examples describes the quotient $C^*_r(P) / \spkl{I}$ as a full corner in a crossed product which admits a finite length independent resolution. Moreover $G$ acts freely on $\cJ_{P \subseteq G}\reg$, $G\setminus\cJ_{P \subseteq G}^\times$ is a singleton, and $\sR$ satisfies conditions (A)-(C) of \cite[\S~5]{L-N} with $i\#j=j$ for all $i\neq j$. We can now use Proposition~\ref{k-theory-free-action-prop} to describe the K-theory of $C^*_r(P) / \spkl{I}$ for $1\leq n\leq 3$. First we see that the matrices $M_i:\Zz_0[G\setminus\cJ_{P \subseteq G}]\to\Zz_0[G\setminus\cJ_{P \subseteq G}]$ defined in \cite[\S~5]{L-N} are given by $[X]\mapsto \sum_{Y\in\cR_i(X)}[Y]$. Since $\Zz_0[G\setminus\cJ_{P \subseteq G}]=\Zz$, we then get $M_i x=p_i x$ for each $x\in\Zz$. As noted in Remark~\ref{chain-complex-remark} we can use the chain complex $\wt{C}$ defined in \cite[\S~5]{L-N} for homology computations. We get for $n=1$ ($p\defeq p_1$),
$$
  C=\wt{C}= \left(
  0 \to \Zz \overset{(1-p)}{\lori} \Zz \to 0
	\right)
$$
and so by Proposition~\ref{k-theory-free-action-prop} and the following remark,
\bglnoz
  K_0(C^*_r(P) / \spkl{I}) &\cong& \Zz/ (1-p)\Zz,\\
	K_1(C^*_r(P) / \spkl{I}) &=& 0.
\eglnoz
For $n=2$ we get
$$
  \wt{C}= \left(
  0 \to \Zz \overset{d_2}{\lori}\Zz\oplus\Zz \overset{d_1}{\lori} \Zz \to 0
	\right)
$$
with
$$
d_2=\begin{bmatrix}p_2-1\\1-p_1\end{bmatrix}, \, d_1= \begin{bmatrix}1-p_1 & 1-p_2\end{bmatrix}.
$$
If we let $g=\gcd(p_1-1,p_2-1)$ this gives us $H_2(\wt{C})=0$, $H_1(\wt{C})=H_0(\wt{C})=\Zz/g\Zz$, so
\bglnoz
  K_0(C^*_r(P) / \spkl{I}) \cong \Zz/g\Zz,\\
	K_1(C^*_r(P) / \spkl{I}) \cong \Zz/g\Zz.
\eglnoz
Moving on to the case $n=3$ we get
$$
  \wt{C}= \left(
  0 \to \Zz \overset{d_3}{\lori} \Zz\oplus\Zz\oplus\Zz \overset{d_2}{\lori} \Zz\oplus\Zz\oplus\Zz \overset{d_1}{\lori} \Zz \to 0
	\right)
$$
with
$$
d_3=\begin{bmatrix}1-p_3\\p_2-1\\1-p_1\end{bmatrix}, \, d_2= \begin{bmatrix}p_2-1 & p_3-1 & 0\\1-p_1 & 0 &p_3-1\\ 0 & 1-p_1 & 1-p_2\end{bmatrix}, \, d_1=\begin{bmatrix}1-p_1 & 1-p_2 & 1-p_3\end{bmatrix}.
$$
Let $g=\gcd(p_1-1,p_2-1,p_3-1)$. Then $H_3(\wt{C})=0$, $H_2(\wt{C})=H_0(\wt{C})=\Zz/g\Zz$ and $H_1(\wt{C})=\Zz/g\Zz\oplus\Zz/g\Zz$. So
$$
  K_1(C^*_r(P) / \spkl{I}) \cong \Zz/g\Zz\oplus\Zz/g\Zz
$$
and there is an extension
$$
0 \to \Zz/g\Zz \to K_0(C^*_r(P) / \spkl{I}) \to \Zz/g\Zz \to 0.
$$

\subsection{C*-algebras of semigroups which do not satisfy independence}

We show that our methods allow us to compute K-theory for semigroup C*-algebras in the case where the independence condition is not satisfied. Let us start with a general observation.

Assume that $D$ is a commutative C*-algebra generated by projections. This means that there exists a semilattice $E$ and a surjective homomorphism $\pi$: $C^*_u(E) \to D$. Further assume that for every $e \in E\reg$, we are given a finite set $\sR(e)$ of finite covers of $e$ such that for every $e \in E\reg$ and $\cR \in \sR(e)$, we have $\pi(e) = \pi(\bigvee\cR)$ in $D$.
\blemma
\label{C/I}
Assume that condition~(i) from Theorem~\ref{maintheo} holds for $E$ and $\sR(e)$, $e \in E\reg$. If for every $e \in E\reg$ and $\gekl{e_i}_{i=1}^n \subseteq E$, $\pi(e) = \pi(\bigvee_{i=1}^n e_i)$ in $D$ implies that there exists $\cR \in \sR(e)$ with $\cR \subseteq \gekl{e_i}_{i=1}^n$, then
\bgloz
  \ker(\pi) = \spkl{\menge{e - \bigvee\cR}{e \in E\reg, \, \cR \in \sR(e)}} \triangleleft C^*_u(E).
\egloz
\elemma
\bproof
Write $I \defeq \spkl{\menge{e - \bigvee\cR}{e \in E\reg, \, \cR \in \sR(e)}} \triangleleft C^*_u(E)$. We obviously have $I \subseteq \ker(\pi)$. To show $I = \ker(\pi)$, we show that the homomorphism $C^*_u(E)/I \to D$ induced by $\pi$ is injective. By \cite[Lemma~2.20]{Li1}, we have to show that for all $d$ and $d_1, \dotsc, d_n$ in $E$, $\pi(d) = \pi(\bigvee_{i=1}^n d_i)$ in $D$ implies that $d - \bigvee_{i=1}^n d_i$ lies in $I$. Let us suppose that we are given $d$ and $d_1, \dotsc, d_n$ in $E$ with $\pi(d) = \pi(\bigvee_{i=1}^n d_i)$ in $D$. By assumption, we can find $\cQ \in \sR(d)$ with $\cQ \subseteq \gekl{d_i}_{i=1}^n$. Let us prove that $\bigvee_{i=1}^n d_i - \bigvee\cQ$ lies in $\Zz$-$\lspan (\menge{e - \bigvee\cR}{e \in E\reg, \, \cR \in \sR(e)})$. We proceed inductively on the number of elements in $\gekl{d_i}_{i=1}^n \setminus \cQ$. The base case $\gekl{d_i}_{i=1}^n = \cQ$ is trivial. Now assume that we have $\cQ \subseteq \gekl{d_i}_{i=1}^{n-1}$ and $\bigvee_{i=1}^{n-1} d_i - \bigvee\cQ$ lies in $\Zz$-$\lspan (\menge{e - \bigvee\cR}{e \in E\reg, \, \cR \in \sR(e)})$. This means that $\bigvee_{i=1}^{n-1} d_i - \bigvee\cQ = \sum \lambda_e (e-\bigvee\cR)$ for some (finitely many) integer coefficients $\lambda_e$. We compute
\bglnoz
  \bigvee_{i=1}^n d_i - \bigvee\cQ &=& \bigvee_{i=1}^{n-1} d_i + d_n - d_n \cdot \rukl{\bigvee_{i=1}^{n-1} d_i} - \bigvee\cQ \\ 
  &=& \rukl{\bigvee_{i=1}^{n-1} d_i - \bigvee\cQ} + d_n - d_n \cdot \rukl{\bigvee\cQ + \sum \lambda_e (e-\bigvee\cR)} \\
  &=& \rukl{\bigvee_{i=1}^{n-1} d_i - \bigvee\cQ} + (d_n - d_n \bigvee\cQ) - \sum \lambda_e (d_n e - d_n \bigvee\cR).
\eglnoz
Since $E$ and $\sR(e)$, $e \in E\reg$ satisfy condition~(i) from Theorem~\ref{maintheo}, we know that $d_n - d_n \bigvee\cQ = d_n d - d_n \bigvee\cQ$ and $d_n e - d_n \bigvee\cR$ are either $0$ or of the form $(d_n d) - \bigvee\cQ'$ or $(d_n e) - \bigvee\cQ''$ for some $\cQ' \in \sR(d_n d)$, $\cQ'' \in \sR(d_n e)$. As $\bigvee_{i=1}^{n-1} d_i - \bigvee\cQ$ is in $\Zz$-$\lspan (\menge{e - \bigvee\cR}{e \in E\reg, \, \cR \in \sR_e})$ by induction hypothesis, we are done.

We have shown that $\bigvee_{i=1}^n d_i - \bigvee\cQ$ lies in $I$. Thus also $d - \bigvee_{i=1}^n d_i = d - \bigvee Q - \rukl{\bigvee_{i=1}^n d_i - \bigvee\cQ}$ lies in $I$.
\eproof

Now let us come to concrete examples of semigroups which do not satisfy independence. Consider the ring $R \defeq \Zz[i \sqrt{3}]$. Its quotient field is given by $Q = \Qz[i \sqrt{3}]$. $R$ is not integrally closed in $Q$. Let $\alpha \defeq \halb (1 + i \sqrt{3})$. $\alpha$ is a primitive sixth root of unity. The integral closure of $R$ is given by $\bar{R} \defeq \Zz[\alpha]$. We have $Q = \Qz[\alpha]$. The multiplicative units in $R$ are given $R^* = \gekl{\pm 1}$, whereas the multiplicative units in $\bar{R}$ are given by $\bar{R}^* = \spkl{\alpha}$. A straightforward computation shows that the fractional ideals of $R$ are given by $\menge{yR}{y \in Q\reg} \cup \menge{y \bar{R}}{y \in Q\reg}$. This is explained in \cite[Example~4.2]{Ste}. As in \cite{Li3}, we set $\cI(R \subseteq Q) \defeq \menge{(x_1 \cdot R) \cap \dotso (x_n \cdot R)}{x_i \in Q\reg}$. As explained in \cite{Li3}, every element of $\cI(R \subseteq Q)$ is a fractional ideal. But in our special case, we have $\bar{R} = \Zz[\alpha] = \halb R \cap \tfrac{\alpha}{2} R \in \cI(R \subseteq Q)$. Thus, the set of fractional ideals coincides with $\cI(R \subseteq Q)$. Moreover, note that $(R : \bar{R}) = \menge{x \in Q}{x \bar{R} \subseteq R} = 2 \bar{R}$. It turns out that $\cI(R \subseteq Q)$ is not independent. Indeed, it is straightforward to see the following
\blemma
\label{relations-ideals}
\begin{itemize}
\item[(a)] We have $\bar{R} = R \cup \alpha R \cup \alpha^2 R$.
\item[(b)] We have $R \cap \alpha R = R \cap \alpha^2 R = \alpha R \cap \alpha^2 R = 2 \bar{R}$, and $2 \bar{R}$ is a proper subset of $R$, $\alpha R$ or $\alpha^2 R$.
\item[(c)] If $\bar{R} = \bigcup_{i=1}^n I_i$ for fractional ideals $I_i$ with $I_i \subsetneq \bar{R}$, then we must have $\gekl{R, \alpha R, \alpha^2 R} \subseteq \menge{I_i}{1 \leq i \leq n}$.
\item[(d)] Let $I$ be a fractional ideal. If $I \cap \bar{R} = yR$ for some $y \in Q\reg$, then $yR \in \gekl{I \cap R, I \cap \alpha R, I \cap \alpha^2 R}$. If $I \cap \bar{R} = y \bar{R}$ for some $y \in Q\reg$, then $y \bar{R} \in \gekl{I \cap R, I \cap \alpha R, I \cap \alpha^2 R}$ or $\gekl{I \cap R, I \cap \alpha R, I \cap \alpha^2 R} = \gekl{y R, y \alpha R, y \alpha^2 R}$.
\end{itemize}
\elemma

Let us turn to semigroup C*-algebras. We start with the multiplicative semigroup $R\reg$. The constructible ideals of $R\reg$ are given by $\cJ(R\reg) = \menge{a R\reg}{a \in R\reg} \cup \menge{2c \bar{R}}{c \in \bar{R}} \cup \gekl{\emptyset}$. $R\reg$ is a subsemigroup of the multiplicative group $Q\reg$, and the constructible $R\reg$-ideals in $Q\reg$ are given by $\cJ(R \reg \subseteq Q\reg) = \menge{y R\reg, y \bar{R}\reg}{y \in Q\reg} \cup \gekl{\emptyset}$. $\cJ(R \reg \subseteq Q\reg)$ is a semilattice under intersection ($X Y \defeq X \cap Y$). Let us set for $y \in Q\reg$: $\sR(y R\reg) \defeq \emptyset$ and $\cR(y \bar{R}\reg) \defeq \gekl{y R\reg, y \alpha R\reg, y \alpha^2 R\reg}$, $\sR(y \bar{R}\reg) \defeq \gekl{\cR(y \bar{R}\reg)}$. Using Lemma~\ref{relations-ideals}, it is easy to see that $\cR(y \bar{R}\reg)$ is a finite cover for $y \bar{R}\reg$, and that $Q\reg \curvearrowright \cJ(R \reg \subseteq Q\reg)$ and $\sR(Y)$, $Y \in \cJ(R \reg \subseteq Q\reg)$, satisfy conditions (i) to (iv) of our Theorem~\ref{maintheo} and the assumptions in Lemma~\ref{C/I}. Thus, if we write $E$ for the semilattice $\cJ(R\reg \subseteq Q\reg)$ from above, and if $D$ is the canonical commutative sub-C*-algebra of $\ell^{\infty}(Q\reg)$ corresponding to $\cJ(R\reg \subseteq Q\reg)$ (see \cite[Definition~3.4]{Li2}), then Lemma~\ref{C/I} tells us that $D \cong C^*_u(E) / I$. Here $I$ is the ideal of $C^*_u(E)$ corresponding to our relations $\sR(Y)$, $Y \in \cJ(R \reg \subseteq Q\reg)$. We are now able to compute K-theory for the reduced semigroup C*-algebra $C^*_r(R\reg)$. We denote the projection in $C^*_u(E)$ corresponding to $X \in \cJ(R\reg \subseteq Q\reg)$ by $e_X$. Also, we let $E_1$ be the semilattice $\menge{e_X - \bigvee_{Y \in \cR} e_Y}{X \in \cJ(R\reg \subseteq Q\reg), \, \cR \in \sR(X)} \cup \gekl{0}$. $E_1$ is a semilattice of projections in $C^*_u(E)$. Theorem~\ref{maintheo} yields that the following sequence is exact (and $Q\reg$-equivariant):
$$ 0 \to C^*_u(E_1) \to C^*_u(E) \to D \to 0. $$
Here, the first homomorphism is induced by the canonical inclusion $E_1 \into C^*_u(E)$, and the second homomorphism is the canonical projection determined by $e_X \ma E_X$. Since the group $Q\reg$ is amenable, hence exact, the following sequence is also exact:
\bgl
\label{E1ED}
  0 \to C^*_u(E_1) \rtimes_r Q\reg \overset{\iota}{\lori} C^*_u(E) \rtimes_r Q\reg \overset{\pi}{\lori} D \rtimes_r Q\reg \to 0.
\egl
Here, $\iota$ and $\pi$ are induced by the homomorphisms from above.

We can now compute K-theory for $D \rtimes_r Q\reg$ using the six-term exact sequence for \eqref{E1ED}. Consider the homomorphisms
\bglnoz
  && C^*(\spkl{\alpha}) \overset{\phi}{\lori} C^*_u(E_1) \rtimes_r Q\reg, \ 
  u_g \ma (e_{\bar{R}\reg} - (e_{R\reg} + e_{\alpha R\reg} + e_{\alpha^2 R\reg} - e_{2 \bar{R}\reg})) u_g\\
  && C^*(R^*) \overset{\psi_{R^*}}{\lori} C^*_u(E) \rtimes_r Q\reg, \ 
  u_g \ma e_{R\reg} u_g\\
  && C^*(\spkl{\alpha}) \overset{\psi_{\spkl{\alpha}}}{\lori} C^*_u(E) \rtimes_r Q\reg, \ 
  u_g \ma e_{\bar{R}\reg} u_g.
\eglnoz
By \cite[Corollary~3.14]{C-E-L2}, $\phi$ induces an isomorphism in K-theory, and also $(\psi_{R^*})_* + (\psi_{\spkl{\alpha}})_*$: $K_*(C^*(R^*)) \oplus K_*(C^*(\spkl{\alpha})) \to K_*(C^*_u(E) \rtimes_r Q\reg)$ is an isomorphism.

Let $\res_{\spkl{\alpha}}^{R^*}$: $K_*(C^*(\spkl{\alpha})) \to K_*(C^*(R^*))$ and $\ind_{R^*}^{\spkl{\alpha}}$: $K_*(C^*(R^*)) \to K_*(C^*(\spkl{\alpha}))$ be the canonical restriction and induction maps. As a direct computation shows, we have $((\psi_{R^*})_* + (\psi_{\spkl{\alpha}})_*)^{-1} \circ \iota_* \circ \phi_* = ( - \res_{\spkl{\alpha}}^{R^*}, \ind_{R^*}^{\spkl{\alpha}} \circ \res_{\spkl{\alpha}}^{R^*})$ as homomorphisms $K_0(C^*(\spkl{\alpha})) \to K_0(C^*(R^*)) \oplus K_0(C^*(\spkl{\alpha}))$. Further computations show that on the whole, we have
\bglnoz
  && K_0(C^*_r(R\reg)) \cong K_0(D \rtimes_r Q\reg) 
  \cong \coker( - \res_{\spkl{\alpha}}^{R^*}, \ind_{R^*}^{\spkl{\alpha}} \circ \res_{\spkl{\alpha}}^{R^*}) \cong \Zz^8 / \Zz^2 \cong \Zz^6\\
  && K_1(C^*_r(R\reg)) \cong K_1(D \rtimes_r Q\reg) \cong \ker( - \res_{\spkl{\alpha}}^{R^*}, \ind_{R^*}^{\spkl{\alpha}} \circ \res_{\spkl{\alpha}}^{R^*}) \cong \Zz^4.
\eglnoz

Let us now discuss the right reduced semigroup C*-algebra of $R \rtimes R\reg$. The constructible left ideals of $R \rtimes R\reg$ are given by $\cJ_{\rho}(R \rtimes R\reg) = \menge{R \times X}{X \in \cJ(R\reg)}$. $R \rtimes R\reg$ is a subsemigroup of the $ax+b$-group $Q \rtimes Q\reg$, and the constructible left $R \rtimes R\reg$-ideals in $Q \rtimes Q\reg$ are given by
$$\cJ_{\rho}(R \rtimes R\reg \subseteq Q \rtimes Q\reg) = \menge{X \cdot g}{X \in \cJ_{\rho}(R \rtimes R\reg), g \in Q \rtimes Q\reg} \cup \gekl{\emptyset}.$$
$\cJ_{\rho}(R \rtimes R\reg \subseteq Q \rtimes Q\reg)$ is a semilattice under intersection ($X Y \defeq X \cap Y$). Let us set for $g \in Q \rtimes Q\reg$: $\sR((R \times R\reg) \cdot g) \defeq \emptyset$ and $\sR((R \times 2 \bar{R}\reg) \cdot g) \defeq \gekl{\cR({(R \times 2 \bar{R}\reg) \cdot g})}$, where $\cR((R \times 2 \bar{R}\reg) \cdot g) \defeq \gekl{(R \times 2 R\reg) \cdot g, (R \times 2 \alpha R\reg) \cdot g, (R \times 2 \alpha^2 R\reg) \cdot g}$. Using Lemma~\ref{relations-ideals}, it is easy to see that $\cR((R \times 2 \bar{R}\reg) \cdot g)$ is a finite cover for $(R \times 2 \bar{R}\reg) \cdot g$, and that $Q \rtimes Q\reg \curvearrowright \cJ_{\rho}(R \rtimes R\reg \subseteq Q \rtimes Q\reg)$ and $\sR_Y$, $Y \in \cJ_{\rho}(R \rtimes R\reg \subseteq Q \rtimes Q\reg)$, satisfy conditions (i) to (iv) of our Theorem~\ref{maintheo} and the assumptions in Lemma~\ref{C/I}. Hence, writing $E$ for the semilattice $\cJ_{\rho}(R \rtimes R\reg \subseteq Q \rtimes Q\reg)$ and $D$ for the commutative C*-algebra corresponding to $\cJ_{\rho}(R \rtimes R\reg \subseteq Q \rtimes Q\reg)$ as above, Lemma~\ref{C/I} tells us that $D \cong C^*_u(E) / I$. Here $I$ is the ideal of $C^*_u(E)$ corresponding to our relations. Again, this allows us to compute K-theory for the right reduced semigroup C*-algebra $C^*_{\rho}(R \rtimes R\reg)$. We let $E_1$ be the semilattice $\menge{e_X - \bigvee_{Y\in \cR}e_Y}{X \in \cJ_{\rho}(R \rtimes R\reg \subseteq Q \rtimes Q\reg), \, \cR \in \sR(X)} \cup \gekl{0} \subseteq \Proj(C^*_u(E))$. The same argument as for the multiplicative semigroup $R\reg$ yields that the following sequence is exact:
$$ 0 \to C^*_u(E_1) \rtimes_r (Q \rtimes Q\reg) \overset{\iota}{\lori} C^*_u(E) \rtimes_r (Q \rtimes Q\reg) \overset{\pi}{\lori} D \rtimes_r (Q \rtimes Q\reg) \to 0. $$
Here, $\iota$ and $\pi$ are the canonical homomorphisms. We can now compute K-theory for $D \rtimes_r (Q \rtimes Q\reg)$ using the six-term exact sequence for this short exact sequence. Consider the homomorphisms
\bglnoz
  && C^*(2 \bar{R} \rtimes \spkl{\alpha}) \overset{\phi}{\lori} C^*_u(E_1) \rtimes_r (Q \rtimes Q\reg)\\ 
  && u_g \ma (e_{R \times \bar{R}\reg} - (e_{R \times R\reg} + e_{R \times \alpha R\reg} + e_{R \times \alpha^2 R\reg} - e_{R \times 2 \bar{R}\reg})) u_g\\
  && C^*(R^*) \overset{\psi_{R \rtimes R^*}}{\lori} C^*_u(E) \rtimes_r (Q \rtimes Q\reg), \ 
  u_g \ma e_{R \times R\reg} u_g\\
  && C^*(2 \bar{R} \rtimes \spkl{\alpha}) \overset{\psi_{2 \bar{R} \rtimes \spkl{\alpha}}}{\lori} C^*_u(E) \rtimes_r (Q \rtimes Q\reg), \ 
  u_g \ma e_{R \times \bar{R}\reg} u_g.
\eglnoz
By \cite[Corollary~3.14]{C-E-L2}, $\phi$ induces an isomorphism in K-theory, and also
\bglnoz
  (\psi_{R \rtimes R^*})_* &+& (\psi_{2 \bar{R} \rtimes \spkl{\alpha}})_*: \\
  && K_*(C^*(R \rtimes R^*)) \oplus K_*(C^*(2 \bar{R} \rtimes \spkl{\alpha})) \to K_*(C^*_u(E) \rtimes_r (Q \rtimes Q\reg))
\eglnoz
is an isomorphism.

Let $\res_{2 \bar{R} \rtimes \spkl{\alpha}}^{2 \bar{R} \rtimes R^*}$: $K_*(C^*(2 \bar{R} \rtimes \spkl{\alpha})) \to K_*(C^*(2 \bar{R} \rtimes R^*))$, $\ind_{2 \bar{R} \rtimes R^*}^{R \rtimes R^*}$: $K_*(C^*(2 \bar{R} \rtimes R^*)) \to K_*(C^*(R \rtimes R^*))$ and $\ind_{2 \bar{R} \rtimes R^*}^{2 \bar{R} \rtimes \spkl{\alpha}}$: $K_*(C^*(2 \bar{R} \rtimes R^*)) \to K_*(C^*(2 \bar{R} \rtimes \spkl{\alpha}))$ be the canonical restriction and induction maps. Moreover, let $\nu$: $C^*(2 \bar{R} \rtimes \spkl{\alpha}) \to C^*(2 \bar{R} \rtimes \spkl{\alpha})$ be the homomorphism induced by the group homomorphism $2 \bar{R} \rtimes \spkl{\alpha} \to 2 \bar{R} \rtimes \spkl{\alpha}$, $(z,y) \ma (2z,y)$. As a direct computation shows, we have
\bglnoz
  &&((\psi_{R \rtimes R^*})_* + (\psi_{2 \bar{R} \rtimes \spkl{\alpha}})_*)^{-1} \circ \iota_* \circ \phi_* \\
  &=& ( - \ind_{2 \bar{R} \rtimes R^*}^{R \rtimes R^*} \circ \res_{2 \bar{R} \rtimes \spkl{\alpha}}^{2 \bar{R} \rtimes R^*}, \id + \nu_* \circ \ind_{2 \bar{R} \rtimes R^*}^{2 \bar{R} \rtimes \spkl{\alpha}} \circ \res_{2 \bar{R} \rtimes \spkl{\alpha}}^{2 \bar{R} \rtimes R^*} - \nu_*)
\eglnoz
as homomorphisms $K_0(C^*(2 \bar{R} \rtimes \spkl{\alpha})) \to K_0(C^*(R \rtimes R^*)) \oplus K_0(C^*(2 \bar{R} \rtimes \spkl{\alpha}))$. Further computations show that all in all, we have
\bglnoz
  && K_0(C^*_{\rho}(R \rtimes R\reg)) \cong K_0(D \rtimes_r (Q \rtimes Q\reg)) \\
  &\cong& \coker( - \ind_{2 \bar{R} \rtimes R^*}^{R \rtimes R^*} \circ \res_{2 \bar{R} \rtimes \spkl{\alpha}}^{2 \bar{R} \rtimes R^*}, \id + \nu_* \circ \ind_{2 \bar{R} \rtimes R^*}^{2 \bar{R} \rtimes \spkl{\alpha}} \circ \res_{2 \bar{R} \rtimes \spkl{\alpha}}^{2 \bar{R} \rtimes R^*} - \nu_*) \\
  &\cong& \Zz^{16} / \Zz^4 \cong \Zz^{12}, \\
  && K_1(C^*_{\rho}(R \rtimes R\reg)) \cong K_1(D \rtimes_r (Q \rtimes Q\reg)) \\
  &\cong& \ker( - \ind_{2 \bar{R} \rtimes R^*}^{R \rtimes R^*} \circ \res_{2 \bar{R} \rtimes \spkl{\alpha}}^{2 \bar{R} \rtimes R^*}, \id + \nu_* \circ \ind_{2 \bar{R} \rtimes R^*}^{2 \bar{R} \rtimes \spkl{\alpha}} \circ \res_{2 \bar{R} \rtimes \spkl{\alpha}}^{2 \bar{R} \rtimes R^*} - \nu_*) \cong \Zz^6.
\eglnoz

Finally, we discuss the left reduced semigroup C*-algebra of $R \rtimes R\reg$. The constructible right ideals of $R \rtimes R\reg$ are given by $\cJ_{\lambda}(R \rtimes R\reg) = \menge{(r+I) \times I\reg}{I \in \cI(R)}$, where $\cI(R)$ is the set of integral fractional ideals of $R$. $R \rtimes R\reg$ is a subsemigroup of $Q \rtimes Q\reg$, and the constructible right $R \rtimes R\reg$-ideals in $Q \rtimes Q\reg$ are given by $\cJ_{\lambda}(R \rtimes R\reg \subseteq Q \rtimes Q\reg) = \menge{g \cdot X}{g \in Q \rtimes Q\reg, X \in \cJ_{\lambda}(R \rtimes R\reg)} \cup \gekl{\emptyset}$. $\cJ_{\lambda}(R \rtimes R\reg \subseteq Q \rtimes Q\reg)$ is a semilattice under intersection ($X Y \defeq X \cap Y$). Let us set for $g \in Q \rtimes Q\reg$: $\sR(g \cdot (R \times R\reg)) \defeq \emptyset$ and $\sR(g \cdot (\bar{R} \times \bar{R}\reg)) \defeq \gekl{\cR(g \cdot (\bar{R} \times \bar{R}\reg))}$, where
\bgloz
  \cR(g \cdot (\bar{R} \times \bar{R}\reg)) \defeq \left\{
  \begin{aligned}
    & g \cdot (R \times R\reg), g \cdot ((\alpha + R) \times R\reg), \\
    & g \cdot (\alpha R \times \alpha R\reg), g \cdot ((1 + \alpha R) \times \alpha R\reg), \\
    & g \cdot (\alpha^2 R \times \alpha^2 R\reg), g \cdot ((1 + \alpha^2 R) \times \alpha^2 R\reg) 
  \end{aligned}
  \right\}.
\egloz
Again, using Lemma~\ref{relations-ideals}, it is easy to see that $\cR(g \cdot (\bar{R} \times \bar{R}\reg))$ is a finite cover for $g \cdot (\bar{R} \times \bar{R}\reg)$, and that $Q \rtimes Q\reg \curvearrowright \cJ_{\lambda}(R \rtimes R\reg \subseteq Q \rtimes Q\reg)$ and $\sR(Y)$, $Y \in \cJ_{\lambda}(R \rtimes R\reg \subseteq Q \rtimes Q\reg)$, satisfy conditions (i) to (iv) of our Theorem~\ref{maintheo} and the assumptions in Lemma~\ref{C/I}. Hence, writing $E$ for the semilattice $\cJ_{\lambda}(R \rtimes R\reg \subseteq Q \rtimes Q\reg)$ and $D$ for the commutative C*-algebra corresponding to $\cJ_{\lambda}(R \rtimes R\reg \subseteq Q \rtimes Q\reg)$ as above, Lemma~\ref{C/I} tells us that $D \cong C^*_u(E) / I$. Here $I$ is the ideal of $C^*_u(E)$ corresponding to our relations. We let $E_1$ be the semilattice $\menge{e_X - \bigvee_{Y\in\cR}e_Y}{X \in \cJ_{\lambda}(R \rtimes R\reg \subseteq Q \rtimes Q\reg), \, \cR \in \sR(X)} \cup \gekl{0} \subseteq \Proj(C^*_u(E))$. As before, we obtain that the following sequence is exact:
$$ 0 \to C^*_u(E_1) \rtimes_r (Q \rtimes Q\reg) \overset{\iota}{\lori} C^*_u(E) \rtimes_r (Q \rtimes Q\reg) \overset{\pi}{\lori} D \rtimes_r (Q \rtimes Q\reg) \to 0, $$
where $\iota$ and $\pi$ are the canonical homomorphisms. We can now compute K-theory for $D \rtimes_r (Q \rtimes Q\reg)$ using the six-term exact sequence for this short exact sequence.
Let $\ve$ be given by
\bglnoz
  &=& e_{R \times R\reg} + e_{(\alpha + R) \times R\reg}
  + e_{\alpha R \times \alpha R\reg} + e_{(1 + \alpha R) \times \alpha R\reg}
  + e_{\alpha^2 R \times \alpha^2 R\reg} + e_{(1 + \alpha^2 R) \times \alpha^2 R\reg}\\
  &-& 
   (
   e_{2 \bar{R} \times 2 \bar{R}\reg} + e_{(1 + 2 \bar{R}) \times 2 \bar{R}\reg} 
   + e_{(\alpha + 2 \bar{R}) \times 2 \bar{R}\reg} + e_{(1 + \alpha + 2 \bar{R}) \times 2 \bar{R}\reg}
   ).
\eglnoz
Consider the homomorphisms
\bglnoz
  && C^*(2 \bar{R} \rtimes \spkl{\alpha}) \overset{\phi}{\lori} C^*_u(E_1) \rtimes_r (Q \rtimes Q\reg), \, u_g \ma (e_{\bar{R} \times \bar{R}\reg} - \ve) u_g, \\
  && C^*(R^*) \overset{\psi_{R \rtimes R^*}}{\lori} C^*_u(E) \rtimes_r (Q \rtimes Q\reg), \ 
  u_g \ma e_{R \times R\reg} u_g,\\
  && C^*(\bar{R} \rtimes \spkl{\alpha}) \overset{\psi_{\bar{R} \rtimes \spkl{\alpha}}}{\lori} C^*_u(E) \rtimes_r (Q \rtimes Q\reg), \ 
  u_g \ma e_{\bar{R} \times \bar{R}\reg} u_g.
\eglnoz
By \cite[Corollary~3.14]{C-E-L2}, $\phi$ induces an isomorphism in K-theory, and also
$$
  (\psi_{R \rtimes R^*})_* + (\psi_{\bar{R} \rtimes \spkl{\alpha}})_*: \:
  K_*(C^*(R \rtimes R^*)) \oplus K_*(C^*(\bar{R} \rtimes \spkl{\alpha})) \to K_*(C^*_u(E) \rtimes_r (Q \rtimes Q\reg))
$$
is an isomorphism.

Let $\res_{\bar{R} \rtimes \spkl{\alpha}}^{R \rtimes R^*}$: $K_*(C^*(\bar{R} \rtimes \spkl{\alpha})) \to K_*(C^*(R \rtimes R^*))$, $\res_{\bar{R} \rtimes \spkl{\alpha}}^{2 \bar{R} \rtimes R^*}$: $K_*(C^*(\bar{R} \rtimes \spkl{\alpha})) \to K_*(C^*(2 \bar{R} \rtimes R^*))$, $\res_{\bar{R} \rtimes \spkl{\alpha}}^{2 \bar{R} \rtimes \spkl{\alpha}}$: $K_*(C^*(\bar{R} \rtimes \spkl{\alpha})) \to K_*(C^*(2 \bar{R} \rtimes \spkl{\alpha}))$ and $\ind_{\bar{R} \rtimes R^*}^{\bar{R} \rtimes \spkl{\alpha}}$: $K_*(C^*(\bar{R} \rtimes R^*)) \to K_*(C^*(\bar{R} \rtimes \spkl{\alpha}))$ be the canonical restriction and induction maps. Moreover, let $\mu$: $C^*(2 \bar{R} \rtimes R^*) \to C^*(\bar{R} \rtimes R^*)$ be the isomorphism induced by the group isomorphism $2 \bar{R} \rtimes R^* \to \bar{R} \rtimes R^*$, $(z,y) \ma (2^{-1}z,y)$, and let $\mu'$: $C^*(2 \bar{R} \rtimes \spkl{\alpha}) \to C^*(\bar{R} \rtimes \spkl{\alpha})$ be the isomorphism induced by the group isomorphism $2 \bar{R} \rtimes \spkl{\alpha} \to \bar{R} \rtimes \spkl{\alpha}$, $(z,y) \ma (2^{-1}z,y)$. As a direct computation shows, we have
\bglnoz
  &&((\psi_{R \rtimes R^*})_* + (\psi_{\bar{R} \rtimes \spkl{\alpha}})_*)^{-1} \circ \iota_* \circ \phi_* \\
  &=& ( - \res_{\bar{R} \rtimes \spkl{\alpha}}^{R \rtimes R^*}, \id + \ind_{\bar{R} \rtimes R^*}^{\bar{R} \rtimes \spkl{\alpha}} \circ \mu_* \circ \res_{\bar{R} \rtimes \spkl{\alpha}}^{2 \bar{R} \rtimes R^*} - \mu'_* \circ \res_{\bar{R} \rtimes \spkl{\alpha}}^{2 \bar{R} \rtimes \spkl{\alpha}})
\eglnoz
as homomorphisms $K_0(C^*(\bar{R} \rtimes \spkl{\alpha})) \to K_0(C^*(R \rtimes R^*)) \oplus K_0(C^*(\bar{R} \rtimes \spkl{\alpha}))$. Further computations show that on the whole, we have
\bglnoz
  && K_0(C^*_{\lambda}(R \rtimes R\reg)) \cong K_0(D \rtimes_r (Q \rtimes Q\reg)) \\
  &\cong& \coker( - \res_{\bar{R} \rtimes \spkl{\alpha}}^{R \rtimes R^*}, \id + \ind_{\bar{R} \rtimes R^*}^{\bar{R} \rtimes \spkl{\alpha}} \circ \mu_* \circ \res_{\bar{R} \rtimes \spkl{\alpha}}^{2 \bar{R} \rtimes R^*} - \mu'_* \circ \res_{\bar{R} \rtimes \spkl{\alpha}}^{2 \bar{R} \rtimes \spkl{\alpha}}) \\
  &\cong& \Zz^{16} / \Zz^4 \cong \Zz^{12}, \\
  && K_1(C^*_{\lambda}(R \rtimes R\reg)) \cong K_1(D \rtimes_r (Q \rtimes Q\reg)) \\
  &\cong& \ker( - \res_{\bar{R} \rtimes \spkl{\alpha}}^{R \rtimes R^*}, \id + \ind_{\bar{R} \rtimes R^*}^{\bar{R} \rtimes \spkl{\alpha}} \circ \mu_* \circ \res_{\bar{R} \rtimes \spkl{\alpha}}^{2 \bar{R} \rtimes R^*} - \mu'_* \circ \res_{\bar{R} \rtimes \spkl{\alpha}}^{2 \bar{R} \rtimes \spkl{\alpha}}) \cong \Zz^6.
\eglnoz

\bremark
As in \cite[\S~6.4]{C-E-L2}, we see that the K-theories of the left and right reduced semigroup C*-algebras of $R \rtimes R\reg$ coincide.
\eremark

\end{document}